\documentclass[a4paper,12pt]{article}
\usepackage{setspace}
\usepackage{typearea}
\usepackage[numbers]{natbib}
\usepackage{amsmath,amssymb}
\usepackage{amsfonts}
\usepackage[english]{babel}
\usepackage{times}
\usepackage{subfig}
\usepackage[mathscr]{eucal}
\usepackage{url}
\usepackage{empheq}
\usepackage[titletoc,title]{appendix}
\usepackage{algorithmic,algorithm}
\usepackage[dvips]{graphicx}
\usepackage{amsmath,amsfonts,amssymb,amsthm,epsfig,epstopdf,titling,url,array}
\usepackage{authblk}
\theoremstyle{plain}
\newtheorem{thm}{Theorem}[section]
\newtheorem{lem}[thm]{Lemma}

\newtheorem{cor}[thm]{Corollary}

\theoremstyle{definition}

\theoremstyle{remark}
\newtheorem*{rem}{Remark}

\newtheorem*{pr}{Proof}
\newtheorem*{pr main thm1}{Proof of Theorem \ref{main thm1}}
\newtheorem*{pr main thm2}{Proof of Theorem \ref{main thm2}}

\theoremstyle{assumption}
\newtheorem{assump}{Assumption}
\typearea{13}
\makeatletter
\renewcommand{\AB@affilsep}{\quad\protect\Affilfont}

\makeatother

\title{Stochastic dual averaging methods using variance reduction techniques for regularized empirical risk minimization problems}
\author[1]{Tomoya Murata\thanks{Email: murata.t.ab@m.titech.ac.jp}}
\author[1, 2]{Taiji Suzuki\thanks{Email: suzuki.t.ct@m.titech.ac.jp}}
\affil[1]{Department of Mathematical and Computing Sciences,
Graduate School of Information Science and Engineering,
Tokyo Institute of Technology}
\affil[2]{PRESTO, Japan Science and Technology Agency (JST)}
\date{}
\begin{document}

\bibliographystyle{abbrvnat}
\maketitle
\begin{abstract}
We consider a composite convex minimization problem associated with regularized empirical risk minimization, which often arises in machine learning. We propose two new stochastic gradient methods that are based on stochastic dual averaging method with variance reduction. Our methods generate a sparser solution than the existing methods because we do not need to take the average of the history of the solutions. This is favorable in terms of both interpretability and generalization. Moreover, our methods have theoretical support for both a strongly and a non-strongly convex regularizer and achieve the best known convergence rates among existing nonaccelerated stochastic gradient methods.   
\end{abstract}
\section{Introduction}
\label{intro} 
We consider the following composite convex minimization problem:
\begin{equation}\label{problem}
\underset{x \in \mathbb{R}^d}{\mathrm{min}}\ \ \{P(x)\overset{\mathrm{def}}{=}F(x)+R(x)\}, 
\end{equation}
where  $F(x)=\frac{1}{n}\sum_{i=1}^n f_{i}(x)$. Here each $f_{i}:\mathbb{R}^d \to \mathbb{R}$ is an $L_i$-smooth convex function and $R:\mathbb{R}^d \to \mathbb{R}$ is a relatively simple and (possibly) nondifferentiable convex function. Problems of this form
 often arise in machine learning and are known as regularized empirical risk minimization.  \par
A traditional method for solving (\ref{problem}) is the (proximal) gradient descent (GD) method. The GD algorithm  is very simple and intuitive and achieves a linear convergence rate for a strongly convex regularizer. However, in typical machine learning tasks, the number $n$ can be very large, and then the iteration cost of GD can be quite expensive. \par
A popular alternative for solving (\ref{problem}) is the stochastic gradient descent (SGD) method \citep{singer2009efficient, hazan2007logarithmic, shalev2007logarithmic}. Since the iteration cost of SGD is very cheap, SGD is suitable to many machine learning tasks. However, SGD only achieves a  sublinear convergence rate and is ultimately slower than GD. \par
Recently, a number of (first-order) stochastic gradient methods using variance reduction techniques, which utilize the finite sum structure of problem (\ref{problem}), have been proposed \citep{roux2012stochastic, schmidt2013minimizing, johnson2013accelerating, xiao2014proximal, nitanda2014stochastic, defazio2014saga, allen2015univr}. The iteration costs of these methods are the same as that of SGD, and, moreover, they achieve a linear convergence rate for a strongly convex objective. \par
The  stochastic average gradient (SAG) method \citep{roux2012stochastic, schmidt2013minimizing} can be used to treat the special case of  problem (\ref{problem}) with $R=0$. To the best of our knowledge, SAG  is the first variance reduction algorithm that achieves a linear convergence rate for a strongly convex objective. SAGA \citep{defazio2014saga} is a modified SAG algorithm that not only achieves a linear convergence rate for a strongly convex objective but also can handle a nondifferentiable and non-strongly convex regularizer. However, for a non-strongly convex regularizer, SAGA needs to output the average of the whole history of the solutions for a convergence guarantee whereas SAG and SAGA do not for a strongly convex objective. \par
In contrast, the stochastic variance reduced gradient (SVRG) method \citep{johnson2013accelerating, xiao2014proximal} adopts a different variance reduction scheme from SAG and SAGA, and in Acc-SVRG \citep{nitanda2014stochastic}  a momentum scheme is applied to SVRG. These methods do not have theoretical support for a non-strongly convex regularizer but they achieve a linear convergence rate for a strongly convex objective. (SVRG needs to output the average of the generated solutions in the last stage for a convergence guarantee whereas Acc-SVRG does not.) UniVR \citep{allen2015univr} is an extension of SVRG and can handle a non-strongly convex regularizer and achieves an $O\left(n\mathrm{log}\frac{1}{\varepsilon}+\frac{\bar L}{\varepsilon}\right)$ rate (where  the $O$ notation means the order of the necessary number of the gradient evaluations), which is faster than the  $O\left(\frac{n+L_{\mathrm{max}}}{\varepsilon}\right)$ rate of SAGA  for $\mathbb{E}[P(x)-P(x_*)] \leq \varepsilon$, $\bar L = (1/n) \sum_{i=1}^{n}L_i$, and $L_{\mathrm{max}} = {\mathrm{max}}\{L_1, \ldots , L_n\}$. However,  UniVR also needs to output the average of the generated solutions in the last stage for convergence guarantees for both  strongly and  non-strongly convex regularizers.  \par
In summary, the algorithms used in these methods often need to output the average of the history of the solutions as a final solution for convergence guarantees (and, especially, for a non-strongly convex regularizer, all of these methods need to take the average). This requirement is unsatisfactory for a sparsity-inducing regularizer because the average of the previous solutions could be nonsparse. \par
In this paper, we propose two new stochastic gradient methods using variance reduction techniques: the stochastic variance reduced dual averaging (SVRDA) method and the stochastic average dual averaging (SADA) method. Compared to previous stochastic optimization methods, the main advantages of our algorithms are as follows:
\begin{itemize}
\item Nice sparsity recovery performance: Our algorithms do not need to take the average of the history of the solutions whereas the existing ones do. This property often leads to sparser solutions than the existing methods for sparsity-inducing regularizers.  
\item Fast convergence: Our algorithms achieve the best known convergence rates among the existing nonaccelerated stochastic gradient methods for both  strongly and non-strongly convex regularizers. Experimentally, our algorithms show comparable or superior convergence speed to that of the existing methods. 
\end{itemize}
\section{Assumptions and notation}
\label{assump}
We make the following assumptions for our theory:
\begin{assump}\label{assump1}
Each $f_i$ is convex and differentiable, and its gradient is $L_i$-Lipschitz continuous, i.e., 
\begin{align}
||\nabla f_i(x)-\nabla f_i(y)||_2 \leq L_i||x-y||_2\hspace{0.5cm}(\forall x, y \in \mathbb{R}^d). \label{L-smooth}
\end{align}
\end{assump}
Condition (\ref{L-smooth}) is equivalent to the following conditions (see \citep{nesterov2013introductory}): 
$$f_i (y) \leq f_i(x) + \langle y-x, \nabla f_i(x)\rangle + \frac{L_i}{2}||x-y||_2^2\hspace{0.5cm}(\forall x, y \in \mathbb{R}^d) $$
and 
$$f_i(x) +  \langle y-x, \nabla f_i(x)\rangle + \frac{1}{2L_i}||\nabla f_i(x) - \nabla f_i(y)||_2^2 \leq f_i(y)\hspace{0.5cm}(\forall x, y \in \mathbb{R}^d).$$
\begin{assump}\label{assump2}
The regularization function $R$ is $\mu$-strongly convex (and it is possible that $\mu =0$), i.e., 
$$R(y) \geq R(x) + \xi^{T}(y-x) + \frac{\mu}{2}||y-x||_2^2 \hspace{0.5cm}(\forall x, y \in \mathbb{R}^d, \forall \xi \in \partial R(x)), $$
where $\partial R(x)$ denotes the set of the subgradients of $R$ at $x$. 
\end{assump}
Observe that, if the regularization function $R$ is $\mu$-strongly convex, then the objective function $P$ is also $\mu$-strongly convex. It is well known that a strongly convex function with $\mu > 0$ has a unique minimizer.
\begin{assump}\label{assump3}
The regularization function $R$ is relatively simple, which means that the proximal mapping of $R,$
$$\mathrm{prox}_{R}(y)=\underset{x \in \mathbb{R}^d}{\mathrm{argmin}\ } \left\{ \frac{1}{2}||x-y||_{2}^{2} + R(x) \right\},$$
can be efficiently computed. 
\end{assump}
Since the function $(1/2)||x-y||_{2}^{2} + R(x)$ is $1+\mu$-strongly convex, the function $\mathrm{prox}_{R}$ is well defined regardless of the strong convexity of $R$. Note that $R$ is not necessarily differentiable. 
\begin{assump}\label{assump4}
There exists a minimizer $x_*$ of  problem (\ref{problem}).
\end{assump}
In addition, we define $\bar L = (1/n) \sum_{i=1}^n L_i $. 
Moreover, we define the probability distribution $Q$ on the set $\{1, 2,  \ldots , n\}$ by $Q = \{q_i\}_{i \in \{1, 2,  \ldots , n\}} = \left\{\frac{L_i}{n\bar L}\right\}_{i \in \{1, 2,  \ldots , n\}}$. This probability distribution is used to randomly pick up a data point in each iteration. By employing nonuniform distribution, we can improve the convergence as in \citep{xiao2014proximal}. \par
Many regularized empirical risk minimization problems in machine learning satisfy these assumptions. For example, given a set of training examples $(a_1, b_1), (a_2, b_2), \ldots , (a_n, b_n)$, where  $a_i \in \mathbb{R}^d$ and $b_i \in \mathbb{R}$, if we set $f_i(x) = (1/2)(a_{i}^{\top}x-b_i)^2$ and $R(x)=\lambda||x||_{1}$, we get Lasso regression. Then the above assumptions are satisfied with $L_i = ||a_i||_2$, $\mu = 0$, and $\mathrm{prox}_{R}(y)=(\mathrm{sign}(y_j)\mathrm{max} \{|y_j|-\lambda, 0\})_{j=1}^d$. If we set $f_i(x)=\mathrm{log}(1+\mathrm{exp}(-b_{i}a_{i}^{\top} x))$ and $R(x)=\lambda_1||x||_1+(\lambda_2/2)||x||_2^2$,  we get logistic elastic net regression. Then the above assumptions are satisfied with $L_i = ||a_i||_2^2/4$, $\mu = \lambda_2$, and $\mathrm{prox}_{R}(y)=(1/(1+\lambda_2))(\mathrm{sign}(y_j)\mathrm{max} \{|y_j|-\lambda_1, 0\})_{j=1}^d$. \par
\section{Related work and our contribution}
In this section, we comment on the relationships between our methods and several closely related methods. \par
Standard methods for solving problem (\ref{problem}) are the GD method and the dual averaging (DA) method \citep{nesterov2009primal}. These methods take the following update rules: 
\begin{align}
x_t=& \mathrm{prox}_{\frac{1}{\eta} R}\left(x_{t-1}-\frac{1}{\eta} \nabla F(x_{t-1})\right )\text{\hspace{1cm}(GD), } \notag \\
x_t =& \mathrm{prox}_{\frac{1}{\eta} tR}\left(x_0 - \frac{1}{\eta} \sum_{\tau =1}^{t} \nabla F(x_{\tau-1})\right) \text{\hspace{0.4cm}(DA), } \notag
\end{align}
where $x_0$ is an initial vector and $1/\eta$ is a constant step size. GD and DA   achieve linear convergence rates for a strongly convex regularizer (where, for DA, we need to borrow a multistage scheme as in \citep{chen2012optimal}). However, when the number of data $n$ is very large,  these methods can be quite expensive because they require $O(nd)$ computation for each update. \par Effective alternatives are the SGD method \citep{singer2009efficient, hazan2007logarithmic, shalev2007logarithmic} and the regularized dual averaging (RDA) method \citep{xiao2009dual}. These methods randomly draw $i$ in $\{1,2,\ldots , n\}$ and use $\nabla f_i$ as an estimator of the full gradient $\nabla F$ in each iteration:
\begin{align}
x_t=& \mathrm{prox}_{\frac{1}{\eta_t} R}\left(x_{t-1}-\frac{1}{\eta_t} \nabla f_{i_t}(x_{t-1}) \right)\text{\hspace{1cm}(SGD), } \notag \\
x_t =& \mathrm{prox}_{\frac{1}{\eta_t} tR}\left(x_0 - \frac{1}{\eta_t} \sum_{\tau =1}^{t} \nabla f_{i_{\tau}}(x_{\tau-1})\right) \text{\hspace{0.32cm}(RDA), } \notag
\end{align}
where $1/\eta_t$ is a decreasing step size. These methods only require $O(d)$ computation for each iteration and are suitable for large-scale problems in machine learning. However, though $\nabla f_i$ is an unbiased estimator of $\nabla F$,  it generally has a large variance, which causes  slow convergence. As a result, these methods only achieve sublinear convergence rates even when the regularizer is strongly convex. One of  simple solutions of this problem is to use a mini-batch strategy \citep{cotter2011better, dekel2012optimal}. However, a mini-batch strategy still gives sublinear convergence. \par
In recent years, a number of (first-order) stochastic gradient methods using variance reduction techniques, which utilize the finite sum structure of problem (\ref{problem}), have been proposed \citep{roux2012stochastic, schmidt2013minimizing, johnson2013accelerating, xiao2014proximal, nitanda2014stochastic, defazio2014saga, allen2015univr}. These methods apply a variance reduction technique to SGD. For example, SVRG \citep{johnson2013accelerating, xiao2014proximal}  takes the following update rules: 
\begin{align}
&\widetilde{x} = \widetilde{x}_{s-1} \notag \\
&\text{for } t=1 \text{ to } m \notag \\
&\hspace{0.5cm} \text{Draw } i_t \text{ randomly from } \{1,2,\ldots ,n\} \notag \\
&\hspace{0.5cm} v_t = \nabla f_{i_t}(x_{t-1}) - \nabla f_{i_t}(\widetilde{x}) + \nabla F( \widetilde{x})  \notag \\
&\hspace{0.5cm} x_t = \mathrm{prox}_{\frac{1}{\eta} R}\left(x_{t-1}-\frac{1}{\eta} v_t \right) \notag \\ 
&\widetilde{x}_{s} = \frac{1}{m}\sum_{t=1}^m x_t.  \notag
\end{align}
$v_t$ is an unbiased estimator of $\nabla F(x_{t-1})$ and one can show that  its variance is ``reduced'':
$$\mathrm{E}||v_t - \nabla F(x_{t-1})||^2 \leq 4\bar L [P(x_{t-1}) - P(x_*) + P(\widetilde{x}) - P(x_*)].$$
This means that the variance of the estimator $v_t$ converges to zero as $x_t$ and $\widetilde{x}$ to $x_*$.  In this sense, $v_t$ is a better estimator of $\nabla F(x_{t-1})$ than the simple estimator $\nabla f_{i_t} (x_{t-1})$. 
Indeed, these methods achieve linear convergence rates for a strongly convex regularizer. \par
However, these methods often need to take the average of the previous solutions for convergence guarantee. For example, SVRG and UniVR \citep{allen2015univr} require  taking the average of the history of the solutions in the last stage. SAGA \citep{defazio2014saga} also requires taking the average of all previous solutions for a non-strongly convex regularizer, though it does not for a strongly convex regularizer. For a sparsity-inducing regularizer, this requirement is unsatisfactory because taking the average could cause a nonsparse solution even though the optimal solution is sparse. \par
In contrast, our proposed methods have theoretical convergence guarantees without taking the average of the previous solutions for both  strongly and  non-strongly convex regularizers. The basic idea of our methods is simple: We apply a variance reduction technique to RDA rather than to SGD.  For example, using an analogy to SVRG, we naturally get the following algorithm: 
\begin{align}
&\widetilde{x} = \widetilde{x}_{s-1} \notag \\
&\text{for } t=1 \text{ to } m \notag \\
&\hspace{0.5cm} \text{Draw } i_t \text{ randomly from } \{1,2,\ldots ,n\}  \notag \\
&\hspace{0.5cm} v_t = \nabla f_{i_t}(x_{t-1}) - \nabla f_{i_t}(\widetilde{x}) + \nabla F( \widetilde{x})  \notag \\
&\hspace{0.5cm} x_t =  \mathrm{prox}_{\frac{1}{\eta} tR}\left(x_0 - \frac{1}{\eta} \sum_{\tau =1}^{t}  v_{\tau}\right) \notag \\ 
&\widetilde{x}_{s} = \frac{1}{m}\sum_{t=1}^m x_t.  \notag
\end{align}
However, this algorithm is not sufficient because the final solution has to be the average of the previous solutions for convergence guarantees (a situation that is similar to RDA). Hence we borrow a momentum scheme and an additional SGD step. (For more detail, see Section \ref{algorithm}.) Then the algorithm does not need to take the average of the previous solutions for convergence guarantees even when the regularizer is non-strongly convex. We call this algorithm SVRDA. Similarly, we can apply the dual averaging scheme to SAGA and we call this algorithm SADA. \par
 Comparisons of the properties of these methods are summarized in Table \ref{table}. 
``Gradient complexity'' indicates  the order of the number of the necessary gradient evaluations for $\mathbb{E}[P(x)-P(x_*)] \leq \varepsilon$ (or $\mathbb{E}||x-x_*||_2^2 \leq \varepsilon$). 
``Final output'' indicates whether the (theoretically guaranteed) final solution  is generated from the (weighted) average of previous iterates (Avg) or from the proximal mapping (Prox). For sparsity-inducing regularizers, the solution generated from the proximal mapping is often sparser than the averaged solution. 
As we can see from Table \ref{table}, the proposed SVRDA and SADA both possess  good properties in comparison with state-of-the-art stochastic gradient methods.
\begin{table}
\begin{center}\scriptsize{
\begin{tabular}{|c|c|c|c|c|c|}
\hline
& \multicolumn{2}{c|}{Strongly convex} & \multicolumn{2}{c|}{Non-strongly convex} &  \\ \cline{2-5}
&Gradient complexity&Final output&Gradient complexity&Final output&Memory cost \\ \hline
SAG \citep{roux2012stochastic, schmidt2013minimizing}&$O\left(\left(n+\frac{L_{\mathrm{max}}}{\mu}\right)\mathrm{log} \frac{1+L_{\mathrm{max}}/n}{\varepsilon}\right)$&Prox&$O\left(\frac{n+L_{\mathrm{max}}}{\varepsilon}\right)$&Avg&$O(nd)$ \\ \hline
SVRG \citep{johnson2013accelerating, xiao2014proximal}&$O\left(\left(n+\frac{\bar L}{\mu}\right)\mathrm{log}\frac{1}{\varepsilon}\right)$&Avg&\multicolumn{2}{c|}{No  direct analysis}&$O(d)$ \\ \hline
Acc-SVRG \citep{nitanda2014stochastic}  &$O\left(\left(n+\mathrm{min}\left\{ \frac{\bar L}{\mu}, n\sqrt{\frac{\bar L}{\mu}}\right\} \right)\mathrm{log}\frac{1}{\varepsilon} \right)$&Prox&\multicolumn{2}{c|}{No  direct analysis}&$O(d)$ \\ \hline
SAGA \citep{defazio2014saga}&$O\left(\left(n+\frac{L_{\mathrm{max}}}{\mu}\right)\mathrm{log}\frac{1+n/L_{\mathrm{max}}}{\varepsilon}\right)$&Prox&$O\left(\frac{n+L_{\mathrm{max}}}{\varepsilon}\right)$&Avg&$O(nd)$ \\ \hline
UniVR \citep{allen2015univr}&$O\left(\left(n+\frac{\bar L}{\mu}\right)\mathrm{log}\frac{1}{\varepsilon}\right)$&Avg&$O\left(n\mathrm{log}\frac{1}{\varepsilon}+\frac{\bar L}{\varepsilon}\right)$&Avg&$O(d)$ \\ \hline 
\bf{SVRDA}&$O\left(\left(n+\frac{\bar L}{\mu}\right)\mathrm{log}\frac{1}{\varepsilon}\right)$&Prox&$O\left(n\mathrm{log}\frac{1}{\varepsilon}+\frac{\bar L}{\varepsilon}\right)$&Prox&$O(d)$ \\ \hline
\bf{SADA}&$O\left(\left(n+\frac{L_{\mathrm{max}}}{\mu}\right)\mathrm{log}\frac{1}{\varepsilon}\right)$&Prox&$O\left(n\mathrm{log}\frac{1}{\varepsilon}+\frac{ L_{\mathrm{max}}}{\varepsilon}\right)$&Prox&$O(nd)$ \\ \hline
\end{tabular}}
\end{center}
\caption{Summary of different stochastic gradient methods that use variance reduction techniques}\label{table}
\end{table}
\section{Algorithm description}
\label{algorithm}
In this section, we illustrate the proposed methods. 
\subsection{The SVRDA method}
We provide details  of the SVRDA method in Algorithm \ref{svrda}. The SVRDA method adopts a multistage scheme. Step (\ref{svrg step}) generates  a variance reduced estimator of the full gradient with nonuniform sampling and is the same as  SVRG \citep{johnson2013accelerating, xiao2014proximal}. Update rules (\ref{dual averaging step}) and (\ref{gradient descent step}) are the dual averaging update and the gradient descent update, respectively. The SVRDA method combines these two update rules. This idea is similar to the ORDA method \citep{chen2012optimal}. As in (\ref{itr_num}), for a non-strongly convex regularizer, we have to exponentially increase the iteration number in each inner loop whereas we can use a common fixed iteration number in each inner loop for a strongly convex regularizer. Note that the computational cost of each iteration in the inner loop of the SVRDA method is $O(d)$ rather than $O(nd)$. Also note that SVRDA outputs the solution generated from the proximal mapping rather than the average of previous iterates. 
For a strongly convex regularizer, SVRDA can output both $\widetilde{x}_S$ (the gradient descent step's output) and $\widetilde{v}_S$ (the dual averaging step's output) as a final solution. This is because the convergence of $\mathbb{E}||\widetilde{v}_s-x_*||_2^2$ is guaranteed with a linear convergence rate whereas the theoretical convergence of $\mathbb{E}[P(\widetilde{v}_s)-P(x_*)]$  is not guaranteed. Outputting $\widetilde{v}_S$ experimentally leads to  better sparsity recovery performance than outputting $\widetilde{x}_S$ (see Section \ref{experiments}). 
\begin{algorithm}[t]
\caption{SVRDA}
\label{svrda}
\algsetup{indent = 2em}
\begin{algorithmic}
\REQUIRE $\widetilde{x}_{0} \in \mathbb{R}^d$, $\eta > 0$, $m_1 \in  \mathbb{N}$, $S \in \mathbb{N}$.
\STATE  $\widetilde{v}_{0} = \widetilde{x}_{0}$
\STATE$\alpha = \begin{cases} \frac{1}{4} & (\mu >0) \\ 0 &  (\mu=0)\end{cases}$
\FOR {$s=1$ to $S$}
\STATE $x_0=\widetilde{x}_{s-1}$, $v_0=(1-\alpha)\widetilde{v}_{s-1}+\alpha\widetilde{x}_{s-1}, u_0 = v_0$,  $\bar{g}_0 = 0$
\STATE\begin{equation}m_s =\begin{cases} m_1 &(\mu>0)\\ 2^{s-1} m_1 &(\mu=0) \end{cases} \label{itr_num}\end{equation} 
\hspace{\algorithmicindent}\FOR {$t=1$ to $m_{s}$}
\STATE pick $i_t \in \{1,2,\ldots ,n\}$ randomly according to $Q$
\vspace{-0.6cm}\STATE \begin{align} g_{t}&=(\nabla f_{i_t}(u_{t-1})-\nabla f_{i_t}(x_0))/nq_{{i}_t}+\nabla F(x_0) \hspace{2cm} \label{svrg step}\\
\bar{g}_t &= \left(1-\frac{1}{t}\right)\bar{g}_{t-1} + \frac{1}{t}g_t  \notag \\ 
v_{t}&=\underset{x \in \mathbb{R}^d}{\mathrm{argmin}\ } \left\{ \langle \bar g_t, x \rangle + R(x) + \frac{\eta}{2t}||x-v_0||_2^2 \right\} \notag \\
&=\mathrm{prox}_{\frac{1}{\eta}tR}\left(v_0-\frac{1}{\eta}t\bar{g}_t\right)\label{dual averaging step}\\
x_{t}&= \underset{x \in \mathbb{R}^d}{\mathrm{argmin}\ } \left\{ \langle g_t, x \rangle + R(x) + \frac{\eta t}{2}||x-u_{t-1}||_2^2 \right\} \notag \\
&=\mathrm{prox}_{\frac{1}{\eta t }R}\left(u_{t-1}-\frac{1}{\eta t}g_t\right) \label{gradient descent step}\\
u_t &= \left(1-\frac{1}{t+1}\right)x_t + \frac{1}{t+1}v_t \notag
\end{align}
\vspace{-0.6cm}
\ENDFOR
\STATE $\widetilde{x}_s = x_{m_s}$, $\widetilde{v}_s = v_{m_s}$
\ENDFOR 
\ENSURE $\widetilde{x}_S$ or $\widetilde{v}_S$\ ($\mu>0$), $\widetilde{x}_S$\ ($\mu=0$).
\end{algorithmic}
\end{algorithm}
\subsection{The SADA method}
We provide  details of the SADA method in Algorithm \ref{sada}. The algorithm is similar to  SVRDA  (Algorithm \ref{svrda}). The main difference from the SVRDA method is the update rule (\ref{saga step}). This step reduces the variance of the approximation of the full gradient using a SAGA \citep{defazio2014saga} type variance reduction technique rather than SVRG. Note that SADA is a multistage algorithm like SVRG and SVRDA whereas SAGA is a single-stage algorithm. To the best of our knowledge, there exists no single-stage dual averaging algorithm that achieves a linear convergence rate for a strongly convex regularizer. This is probably because of the limitations of the single-stage dual averaging algorithms. Also note that we adopt uniform sampling for SADA. Schmidt et al. \citep{schmidt2015non} have considered a nonuniform sampling scheme for SAGA on the special setting $R=0$ in (\ref{problem}), but their methods require two  gradient evaluations in one iteration and it is not satisfactory. For this reason, we do not adopt nonuniform sampling schemes for SADA in this paper. SADA has theoretically similar properties to SVRDA except for the difference of the sampling scheme, and experimentally SADA sometimes outperforms SVRDA (see Section \ref{experiments}).  
\begin{algorithm}[t]
\caption{SADA}
\label{sada}
\algsetup{indent = 2em}
\begin{algorithmic}
\REQUIRE $\widetilde{x}_{0} \in \mathbb{R}^d$, $\eta>0$, $m_1 \in  \mathbb{N}$, $S \in \mathbb{N}$
\STATE  $\widetilde{v}_{0} = \widetilde{x}_{0} $
\STATE $\alpha = \begin{cases} \frac{1}{4} & (\mu >0) \\ 0 &  (\mu=0)\end{cases}$
\FOR {$s=1$ to $S$}
\STATE $x_0=\widetilde{x}_{s-1}$, $v_0 = (1-\alpha)\widetilde{v}_{s-1}+\alpha \widetilde{x}_{s-1}$, $u_0 = v_0$,  $\bar{g}_0 = 0$,  $\phi_i^0 = x_0 \ (i=1,2, \ldots , n)$
\STATE $$m_s =\begin{cases} m_1 &(\mu>0)\\ 2^{s-1} m_1 &(\mu=0)\end{cases}$$
\FOR {$t=1$ to $m_{s}$}
\STATE pick $i_t \in \{1,2,\ldots ,n\}$ uniformly at random
\vspace{-0.6cm} \STATE \begin{align} \phi_{i_t}^t &= u_{t-1}, \phi_i^t = \phi_i^{t-1} (i \neq i_t) \notag \\
g_{t}&=\nabla f_{i_t}(\phi_{i_t}^t)-\nabla f_{i_t}(\phi_{i_t}^{t-1})+\frac{1}{n}\sum_{i=1}^{n}\nabla f_i(\phi_{i}^{t-1})  \hspace{2cm} 
\label{saga step}\\
\bar{g}_t &= \left(1-\frac{1}{t}\right)\bar{g}_{t-1} + \frac{1}{t}g_t \notag \\
v_{t}&=\underset{x \in \mathbb{R}^d}{\mathrm{argmin}\ } \left\{ \langle \bar g_t, x \rangle + R(x) + \frac{\eta}{2t}||x-v_0||_2^2 \right\} \notag \\
&=\mathrm{prox}_{\frac{t}{\eta}R}\left(v_0-\frac{t}{\eta}\bar{g}_t\right) \notag \\
x_{t}&= \underset{x \in \mathbb{R}^d}{\mathrm{argmin}\ } \left\{ \langle g_t, x \rangle + R(x) + \frac{\eta t}{2}||x-u_{t-1}||_2^2 \right\} \notag \\
&=\mathrm{prox}_{\frac{1}{\eta t}R}\left(u_{t-1}-\frac{1}{\eta t}g_t\right)\notag \\
u_t &= \left(1-\frac{1}{t+1}\right)x_t + \frac{1}{t+1}v_t\notag 
\end{align}\vspace{-0.6cm}
\ENDFOR
\STATE $\widetilde{x}_s = x_{m_s}$, $\widetilde{v}_s = v_{m_s}$
\ENDFOR 
\ENSURE $\widetilde{x}_S$ or $\widetilde{v}_S$\ ($\mu>0$), $\widetilde{x}_S$\ ($\mu=0$).
\end{algorithmic}
\end{algorithm}
\clearpage
\section{Convergence analysis}
Now we give a convergence analysis of our algorithms. In this section, all norms $||\cdot||$ mean the $L_2$-norm $||\cdot||_2$. 
\subsection{Convergence analysis of SVRDA}
In this subsection, we give the convergence analysis of SVRDA.
\begin{thm}
\label{main thm1}
Suppose that Assumptions \ref{assump1}, \ref{assump2}, \ref{assump3}, and \ref{assump4} hold (and it is possible that $\mu = 0$). Let $x_0 \in \mathbb{R}^d$, $\eta = 4\bar L$, $m_1 \in \mathbb{N,}$ and $0 \leq \alpha \leq 1$. Then the SVRDA algorithm satisfies
\begin{align}
&\mathbb{E} \left[P(\widetilde{x}_{s})-P(x_{*})\right] + \frac{\eta+m_s\mu}{2m_s} \mathbb{E}||\widetilde{v}_{s} - x_*||^2\notag \\ 
\leq& \frac{1}{2} \mathbb{E} [P(\widetilde{x}_{s-1})-P(x_{*})] + \left(\frac{\alpha\eta }{2m_s}-\frac{\mu}{4}\right)\mathbb{E} ||\widetilde{x}_{s-1} - x_*||^2+\frac{(1-\alpha)\eta }{2m_s}\mathbb{E} ||\widetilde{v}_{s-1} - x_*||^2.\notag 
\end{align}
\end{thm}
\begin{rem}On inequality (\ref{roughbound}) in Appendix \ref{appendix}, we can apply a tighter bound and  $\eta$ can be smaller than $4\bar L$ for satisfying Theorem \ref{main thm1}. This means that we can get a larger step size $\frac{1}{\eta}$ than $\frac{1}{4\bar L}$ and have a theoretically tighter bound. However, practically, if we tune $\eta$, it makes little difference and thus we omit it in this paper. 
\end{rem}
The proof of Theorem \ref{main thm1} is given in Appendix \ref{appendix}. Using this theorem, we derive recursive inequalities relative to $\mathbb{E} [P(\widetilde{x}_{s})-P(x_{*})]$ and $\mathbb{E}||\widetilde{v}_{s} - x_*||^2$. Based on Theorem \ref{main thm1}, we obtain the linear convergence of SVRDA for $\mu > 0$. 
\begin{cor}[for a strongly convex regularizer]
\label{corsc}
Suppose that Assumptions \ref{assump1}, \ref{assump2}, \ref{assump3}, and \ref{assump4} hold. Moreover, assume that $\mu > 0$. Let $x_0\in \mathbb{R}^d$, $\eta = 4\bar L$, $m_1 = \frac{\eta}{2\mu}$, $S\in \mathbb{N,}$ and $\alpha = \frac{1}{4}$.
Then the SVRDA algorithm satisfies
$$\mathbb{E} [P(\widetilde{x}_{S})-P(x_{*})] + \frac{3\mu}{2} \mathbb{E}||\widetilde{v}_{S} - x_*||^2 \leq \frac{1}{2^S}\left[P(\widetilde{x}_{0})-P(x_{*}) + \frac{3\mu}{2} ||\widetilde{x}_{0} - x_*||^2\right]. $$
In addition, the SVRDA algorithm has a gradient complexity of $$O\left(\left(n+\frac{\bar L}{\mu}\right)\mathrm{log } \frac{1}{\varepsilon} \right)$$ for $\mathbb{E}[P(\widetilde{x}_S)-P(x_{*})] \leq \varepsilon$ and 
$$O\left(\left(n+\frac{\bar L}{\mu}\right)\mathrm{log } \frac{1}{\mu \varepsilon} \right)$$
for $\mathbb{E}||\widetilde{v}_S-x_*||^2 \leq \varepsilon$. 
\end{cor}
These gradient complexities are essentially the same as the ones obtained by \citep{roux2012stochastic, schmidt2013minimizing, johnson2013accelerating, xiao2014proximal, nitanda2014stochastic, defazio2014saga, allen2015univr} and are the best known ones among the existing nonaccelerated stochastic gradient methods. Note that the gradient complexity of GD is $O\left(n\frac{\bar L}{\mu}\mathrm{log}\frac{1}{\varepsilon}\right)$ and that of SGD is $O\left(\frac{1}{\mu \varepsilon}\right)$. In a typical empirical risk minimization task, we require that $\varepsilon$ be $O\left(\frac{1}{n}\right)$. Then the gradient complexities of SVRDA, GD, and SGD are $O\left(\left(n+\frac{\bar L}{\mu}\right)\mathrm{log } n \right)$, $O\left(n\frac{\bar L}{\mu}\mathrm{log}n \right)$, and $O\left(\frac{n}{\mu}\right)$, respectively.  Hence, SVRDA significantly improves upon the gradient complexities of GD and SGD for $\mu>0$. 
\begin{pr}
By Theorem \ref{main thm1} and the definitions of $\eta$, $m_s,$ and $\alpha$, we obtain
\begin{align}
&\mathbb{E} [P(\widetilde{x}_{S})-P(x_{*})] + \frac{3\mu}{2} \mathbb{E}||\widetilde{v}_{S} - x_*||^2\notag \\
\leq& \frac{1}{2} \left[\mathbb{E}[P(\widetilde{x}_{S-1})-P(x_{*})] + \frac{3\mu}{2} \mathbb{E}||\widetilde{v}_{S-1} - x_*||^2\right] \notag \\
\leq& \cdots \notag \\
\leq& \frac{1}{2^S} \left[P(\widetilde{x}_{0})-P(x_{*}) + \frac{3\mu}{2} ||\widetilde{v}_{0} - x_*||^2\right] \notag \\
=&\frac{1}{2^S}\left[P(\widetilde{x}_{0})-P(x_{*}) + \frac{3\mu}{2} ||\widetilde{x}_{0} - x_*||^2\right]. \notag 
\end{align}
By this inequality, we can see that the order of the necessary number of outer iterations for $\mathbb{E} [P(\widetilde{x}_{S})-P(x_{*})] \leq \varepsilon$ is $O\left(\mathrm{log}\frac{1}{\varepsilon}\right)$ and  the order of the necessary number of outer iterations for $\mathbb{E} ||\widetilde{v}_{S}-x_{*}||^2 \leq \varepsilon$ is $O\left(\mathrm{log}\frac{1}{\mu \varepsilon}\right)$. Finally,  since SVRDA computes $S$ times the full gradient $\nabla F$ and $2m_1=O\left(\frac{\bar L}{\mu}\right)$ times the gradient $\nabla f_i$, the total gradient complexity is $$O\left(\left(n+\frac{\bar L}{\mu}\right)\mathrm{log } \frac{1}{\varepsilon} \right)$$ for $\mathbb{E}[P(\widetilde{x}_S)-P(x_{*})] \leq \varepsilon$ and 
$$O\left(\left(n+\frac{\bar L}{\mu}\right)\mathrm{log } \frac{1}{\mu \varepsilon} \right)$$
for $\mathbb{E}||\widetilde{v}_S-x_*||^2 \leq \varepsilon$. 
\qed
\end{pr}
Next, we derive the convergence rate for $\mu=0$ from Theorem \ref{main thm1} as follows.
\begin{cor}[for a non-strongly convex regularizer]\label{corsvrganonstrong}
Suppose that Assumptions \ref{assump1}, \ref{assump2}, \ref{assump3}, and \ref{assump4} hold (and it is possible that $\mu = 0$). Let $x_0\in \mathbb{R}^d$, $\eta = 4\bar L$, $m_1$, $S\in \mathbb{N},$ and $\alpha = 0$.
Then the SVRDA algorithm satisfies
$$\mathbb{E} [P(\widetilde{x}_{S})-P(x_{*})]\leq \frac{1}{2^S}\left[P(\widetilde{x}_{0})-P(x_{*}) + \frac{4\bar L }{m_{1}}||\widetilde{x}_{0} - x_*||^2\right].$$
In addition, if $m_1=O(\bar L)$, then 
the SVRDA algorithm has a gradient complexity of $$O\left(n\mathrm{log}\frac{1}{\varepsilon}+ \frac{\bar L }{\varepsilon}\right)$$ for $\mathbb{E}[P(\widetilde{x}_S)-P(x_{*})] \leq \varepsilon$.
\end{cor}
The gradient complexity of SVRDA for a non-strongly convex regularizer is the same as that of UniVR \citep{allen2015univr} and is the best known among the existing stochastic gradient methods. Note that the gradient complexities of GD, SGD, and SAGA \citep{defazio2014saga} are $O\left(\frac{\bar Ln}{\varepsilon}\right)$, $O\left(\frac{1}{\varepsilon^2}\right)$, and $O\left(\frac{n+L_{\mathrm{max}}}{\varepsilon}\right)$, respectively.  In a typical empirical risk minimization task, we require that $\varepsilon$ be $O\left(\frac{1}{n}\right)$. Then the gradient complexities of SVRDA, GD, SGD, and SAGA are $O\left(n\mathrm{log}n+ \bar L n \right)$, $O\left(\bar L n^2\right)$, $O\left(n^2 \right)$, and $O\left(n^2+L_{\mathrm{max}}n \right)$, respectively. Hence, SVRDA significantly improves upon the gradient complexities of GD, SGD, and SAGA for $\mu=0$.
\begin{pr}
By Theorem \ref{main thm1} and the definitions of $\eta$, $m_s$, and $\alpha$, we obtain
\begin{align}
&\mathbb{E} [P(\widetilde{x}_{S})-P(x_{*})] + \frac{\eta}{m_{S+1}} \mathbb{E}||\widetilde{v}_{S} - x_*||^2 \notag \\
=&\mathbb{E} [P(\widetilde{x}_{S})-P(x_{*})] + \frac{\eta}{2m_{S}} \mathbb{E}||\widetilde{v}_{S} - x_*||^2 \notag \\
\leq& \frac{1}{2} \left[\mathbb{E}[P(\widetilde{x}_{S-1})-P(x_{*})] + \frac{\eta}{m_{S}}\mathbb{E}||\widetilde{v}_{S-1} - x_*||^2\right] \notag \\
\leq& \cdots \notag \\
\leq& \frac{1}{2^S} \left[P(\widetilde{x}_{0})-P(x_{*}) + \frac{\eta}{m_{1}}||\widetilde{v}_{0} - x_*||^2\right] \notag \\
=&\frac{1}{2^S}\left[P(\widetilde{x}_{0})-P(x_{*}) + \frac{4\bar L}{ m_{1}}||\widetilde{x}_{0} - x_*||^2\right], \notag 
\end{align}
and therefore
$$\mathbb{E} [P(\widetilde{x}_{S})-P(x_{*})]\leq\frac{1}{2^S} \left[P(\widetilde{x}_{0})-P(x_{*}) + \frac{4\bar L}{ m_{1}}||\widetilde{x}_{0} - x_*||^2\right].$$
Thus the order of the necessary number of outer iterations for $\mathbb{E} [P(\widetilde{x}_{S})-P(x_{*})] \leq \varepsilon$ is $O\left(\mathrm{log} \frac{1}{\varepsilon} \right)$.  
Finally, since SVRDA computes $S$ times the full gradient $\nabla F$ and $O\left(\sum_{s=1}^S 2m_s\right) =  
 O\left(2^Sm_1\right)$ times the gradient $\nabla f_i$, the total gradient complexity is 
$$O\left(n S +\sum_{s=1}^S 2m_s \right) = O\left(n S + 2^Sm_1 \right) = O\left(n\mathrm{log}\frac{1}{\varepsilon} + \frac{\bar L}{\varepsilon} \right).$$ 
\qed
\end{pr}
\subsection{Convergence analysis of SADA}
In this subsection, we give the convergence analysis of SADA.
\begin{thm}\label{main thm2}
Suppose that Assumptions \ref{assump1}, \ref{assump2}, \ref{assump3}, and \ref{assump4} hold (and it is possible that $\mu = 0$). Let $x_0\in \mathbb{R}^d$,  $\eta = 5L_{\mathrm{max}}$, $m_1 \in \mathbb{N,}$ and $0\leq \alpha \leq 1$. Then the SADA algorithm satisfies
\begin{align}
&\mathbb{E} [P(\widetilde{x}_{s})-P(x_{*})] + \frac{\eta+m_s\mu}{2m_s} \mathbb{E}||\widetilde{v}_{s} - x_*||^2 \notag \\
\leq&\frac{1}{2}\mathbb{E} \left[P(\widetilde{x}_{s-1})-P(x_{*})\right] + \left(\frac{\alpha\eta }{2m_s}-\frac{\mu}{4}\right)\mathbb{E}||\widetilde{x}_{s-1} - x_*||^2 + \frac{(1-\alpha)\eta}{2m_s}\mathbb{E}||\widetilde{v}_{s-1} - x_*||^2. \notag
\end{align}
\end{thm}
The proof of Theorem \ref{main thm2} is given in Appendix \ref{appendix2}. Using this theorem, we derive recursive inequalities relative to $\mathbb{E} [P(\widetilde{x}_{s})-P(x_{*})]$ and $\mathbb{E}||\widetilde{v}_{s} - x_*||^2$. Based on Theorem \ref{main thm2}, we obtain the linear convergence of SADA for $\mu > 0$.
\begin{cor}[for a strongly convex regularizer]
\label{corsc2}
Suppose that Assumptions \ref{assump1}, \ref{assump2}, \ref{assump3}, and \ref{assump4} hold. Moreover, assume that $\mu > 0$. Let $x_0\in \mathbb{R}^d$, $\eta = 5L_{{\mathrm{max}}}$, $m_1 = \frac{\eta}{2\mu}$, $S \in \mathbb{N,}$ and $\alpha = \frac{1}{4}$.
Then the SADA algorithm satisfies
$$\mathbb{E} [P(\widetilde{x}_{S})-P(x_{*})] + \frac{3\mu}{2} \mathbb{E}||\widetilde{v}_{S} - x_*||^2 \leq \frac{1}{2^S}\left[P(\widetilde{x}_{0})-P(x_{*} )+ \frac{3\mu}{2}||\widetilde{x}_{0} - x_*||^2\right].$$
In addition, the SADA algorithm has a gradient complexity of $$O\left(\left(n+\frac{L_{\mathrm{max}}}{\mu}\right)\mathrm{log } \frac{1}{\varepsilon} \right)$$ for $\mathbb{E}[P(\widetilde{x}_S)-P(x_{*})] \leq \varepsilon$ and 
$$O\left(\left(n+\frac{L_{\mathrm{max}}}{\mu}\right)\mathrm{log } \frac{1}{\mu \varepsilon} \right)$$
for $\mathbb{E}||\widetilde{v}_S-x_*||^2 \leq \varepsilon$. 
\end{cor}
These gradient complexities are essentially same as the ones obtained by \citep{roux2012stochastic, schmidt2013minimizing, johnson2013accelerating, xiao2014proximal, nitanda2014stochastic, defazio2014saga, allen2015univr} and the ones of SVRDA and are the best known among the existing nonaccelerated stochastic gradient methods. Note that the gradient complexity of GD is $O\left(n\frac{\bar L}{\mu}\mathrm{log}\frac{1}{\varepsilon}\right)$ and that of SGD is $O\left(\frac{1}{\mu \varepsilon}\right)$. In a typical empirical risk minimization task, we require that $\varepsilon$ be $O\left(\frac{1}{n}\right)$. Then the gradient complexities of SADA, GD, and SGD are $O\left(\left(n+\frac{L_{\mathrm{max}}}{\mu}\right)\mathrm{log } n \right)$, $O\left(n\frac{\bar L}{\mu}\mathrm{log}n \right)$, and $O\left(\frac{n}{\mu}\right)$, respectively.  Hence, SADA significantly improves upon the gradient complexities of GD and SGD for $\mu>0$. 
\begin{pr}
The proof of Corollary \ref{corsc} is identical to that of Corollary \ref{corsc} and we omit it. \qed
\end{pr}
\begin{cor}[for non-strongly convex cases]
\label{cornsc2}
Suppose that Assumptions \ref{assump1}, \ref{assump2}, \ref{assump3}, and \ref{assump4} hold (and it is possible that $\mu = 0$). Let $x_0\in \mathbb{R}^d$, $\eta = 5 L_{\mathrm{max}}$, $m_1,S \in \mathbb{N}$, and $\alpha=0$.
Then the SADA algorithm satisfies
$$\mathbb{E}[P(\widetilde{x}_S)-P(x_{*})] \leq \frac{1}{2^{S}} \Bigl[ P(\widetilde{x}_0) - P(x_{*}) + \frac{5L_{\mathrm{max}}}{m_1}||\widetilde{x}_0-x_*||^2   \Bigr].$$ 
In addition,  if $m_1=O(L_{\mathrm{max}})$, then the SADA algorithm has a gradient complexity of $$O\Bigl(n\mathrm{log}\frac{1}{\varepsilon}+ \frac{L_{\mathrm{max}}}{\varepsilon}\Bigr)$$ for $\mathbb{E}[P(\widetilde{x}_S)-P(x_{*})] \leq \varepsilon$.
\end{cor}
The gradient complexity of SADA for a non-strongly convex regularizer is the same as that of UniVR \citep{allen2015univr} and SVRDA and is the best known among the existing stochastic gradient methods. Note that the gradient complexities of GD, SGD, and SAGA \citep{defazio2014saga} are $O\left(\frac{\bar Ln}{\varepsilon}\right)$, $O\left(\frac{1}{\varepsilon^2}\right)$, and $O\left(\frac{n+L_{\mathrm{max}}}{\varepsilon}\right)$, respectively.  In a typical empirical risk minimization task, we require that $\varepsilon$ be $O\left(\frac{1}{n}\right)$. Then the gradient complexities of SVRDA, GD, SGD, and SAGA are $O\left(n\mathrm{log}n+ L_{\mathrm{max}} n \right)$, $O\left(\bar L n^2\right)$, $O\left(n^2 \right)$, and $O\left(n^2+L_{\mathrm{max}}n \right),$ respectively. Hence, SADA significantly improves upon the gradient complexities of GD, SGD, and SAGA for $\mu=0$.
\begin{pr}
The proof is the same as that of Corollary \ref{corsvrganonstrong} and we omit it. \qed
\end{pr}
\section{Numerical experiments}
\label{experiments}
In this section, we provide numerical experiments to demonstrate the performances of SVRDA and SADA. We compare our methods with several state-of-the-art stochastic gradient methods: SVRG \citep{johnson2013accelerating, xiao2014proximal}, SAGA \citep{defazio2014saga}, and UniVR \citep{allen2015univr}. For a fair comparison, we compare all different methods using solutions that are theoretically guaranteed. 
We used nonuniform sampling for SVRG \citep{johnson2013accelerating, xiao2014proximal}, UniVR \citep{allen2015univr}, and SVRDA. (Zhu et al. \citep{allen2015univr} have not considered a nonuniform sampling scheme for UniVR, but because there is theoretical justification of nonuniform sampling for UniVR, we adopted nonuniform sampling for UniVR.) However, we used uniform sampling for SAGA \citep{defazio2014saga} and SADA. (Schmidt et al. \citep{schmidt2015non} considered  nonuniform sampling for SAGA on the special setting $R=0$ in (\ref{problem}), but their algorithm require two  gradient evaluations in one iteration for a gradient complexity of $O\left(\left(n+\frac{\bar L}{\mu}\right)\mathrm{log } \frac{1}{\varepsilon} \right)$, and thus in our experiment we adopted uniform sampling for SAGA and SADA.) \par
In this experiments, we focus on the regularized logistic regression problem for binary classification: Given a set of training examples $(a_1, b_1), (a_2, b_2), \ldots , (a_n, b_n)$, where $a_i \in \mathbb{R}^d$ and $b_i \in \left\{+1,-1\right\}$, we find the optimal classifier $x \in \mathbb{R}^d$ by solving$$\underset{x \in \mathbb{R}^d}{\mathrm{min}}\ \ \frac{1}{n} \sum_{i=1}^n \mathrm{log}(1+\mathrm{exp}(-b_{i}a_{i}^{\top}x)) +\lambda_1||x||_1+\frac{\lambda_2}{2}||x||_2^2,$$
where $\lambda_1$ and $\lambda_2$ are regularization parameters. \par
We used three publicly available data sets in the experiments. Their sizes $n$ and dimensions $d$ are listed in Table \ref{datatable}. Each  continuous feature vector in these data sets has been normalized to  zero mean and unit variance. \par
\begin{table}[H]
\begin{center}
\begin{tabular}{|c|c|c|}\hline
Data sets & $n$ & $d$ \\ \hline
covertype\footnotemark & $581,012$ & $54$ \\ \hline
Reuters-21578\footnotemark &5,964&18,933 \\ \hline
sido0\footnotemark & $12,678$ & $4,932$ \\ \hline
\end{tabular}
\end{center}
\caption{Summary of the data sets used in our numerical experiments}
\label{datatable}
\end{table}
\footnotetext[1]{Available at \url{https://archive.ics.uci.edu/ml/datasets/Covertype}.}
\footnotetext[2]{Available at \url{http://www.cad.zju.edu.cn/home/dengcai/Data/TextData.html}. We converted the 65 class classification task into a binary classification.
}
\footnotetext[3]{Available at \url{http://www.causality.inf.ethz.ch/data/SIDO.html}.}
We performed our experiments on a desktop computer (a Windows 7 64-bit machine with an Intel i7-4790 CPU operating at 3.60 GHz and 8 GB of RAM)  and implemented all algorithms in MATLAB 2015a. \par
\begin{figure}[t]
\subfloat[$\lambda_1=10^{-6}$, $\lambda_2=10^{-6}$]
{\label{covertype_sc}
\begin{minipage}{0.5\hsize}
\centering
\includegraphics[bb=0 0 1500 1000, width=9.5cm]{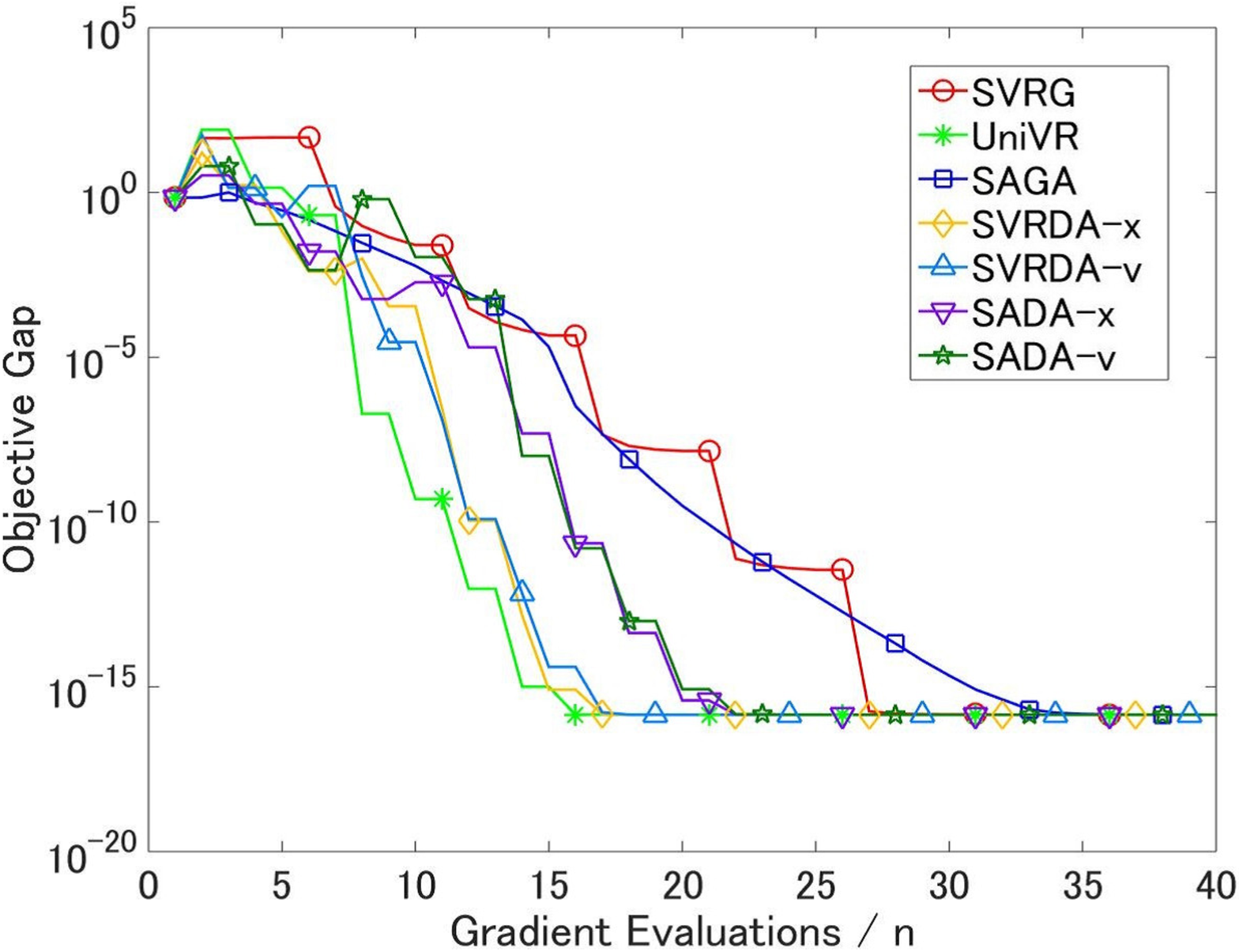}
\end{minipage}
\begin{minipage}{0.5\hsize}
\centering
\includegraphics[bb=0 0 1500 1000, width=9.5cm]{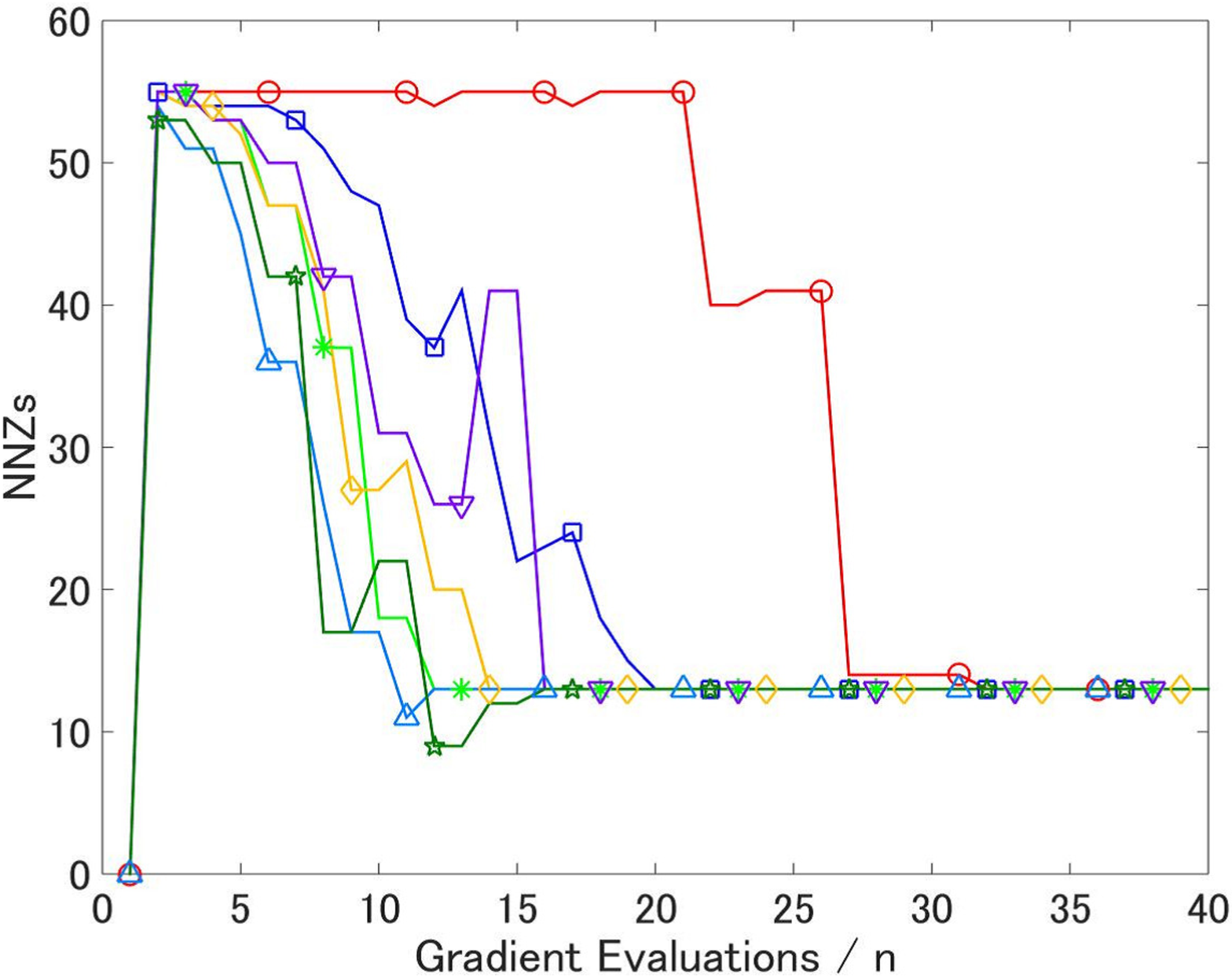}
\end{minipage}
}\\
\subfloat[$\lambda_1=10^{-6}$, $\lambda_2=0$]{\label{covertype_nsc}
\begin{minipage}{0.5\hsize}
\centering
\includegraphics[bb=0 0 1500 1000, width=9.5cm]{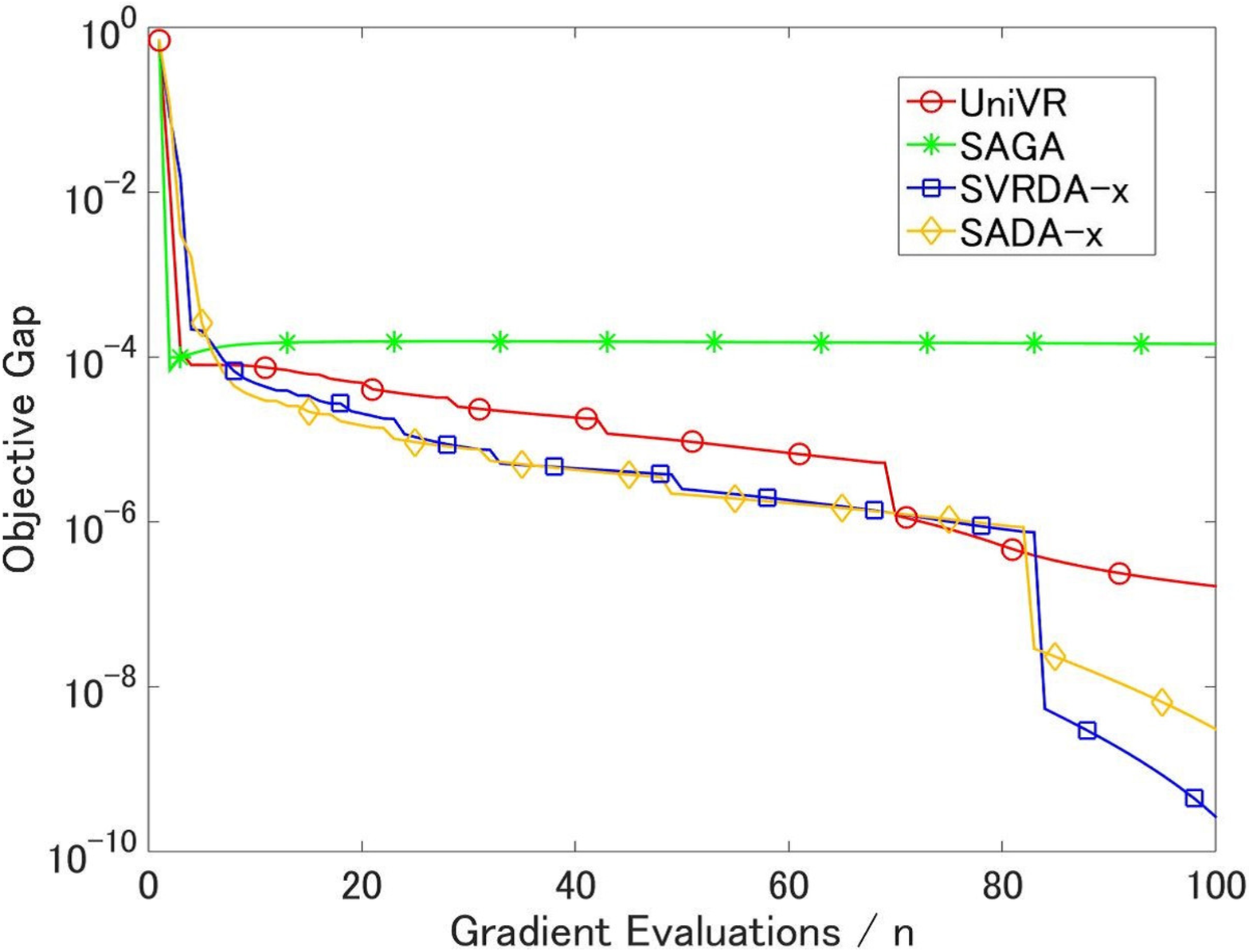}
\end{minipage}
\begin{minipage}{0.5\hsize}
\centering
\includegraphics[bb=0 0 1500 1000, width=9.5cm]{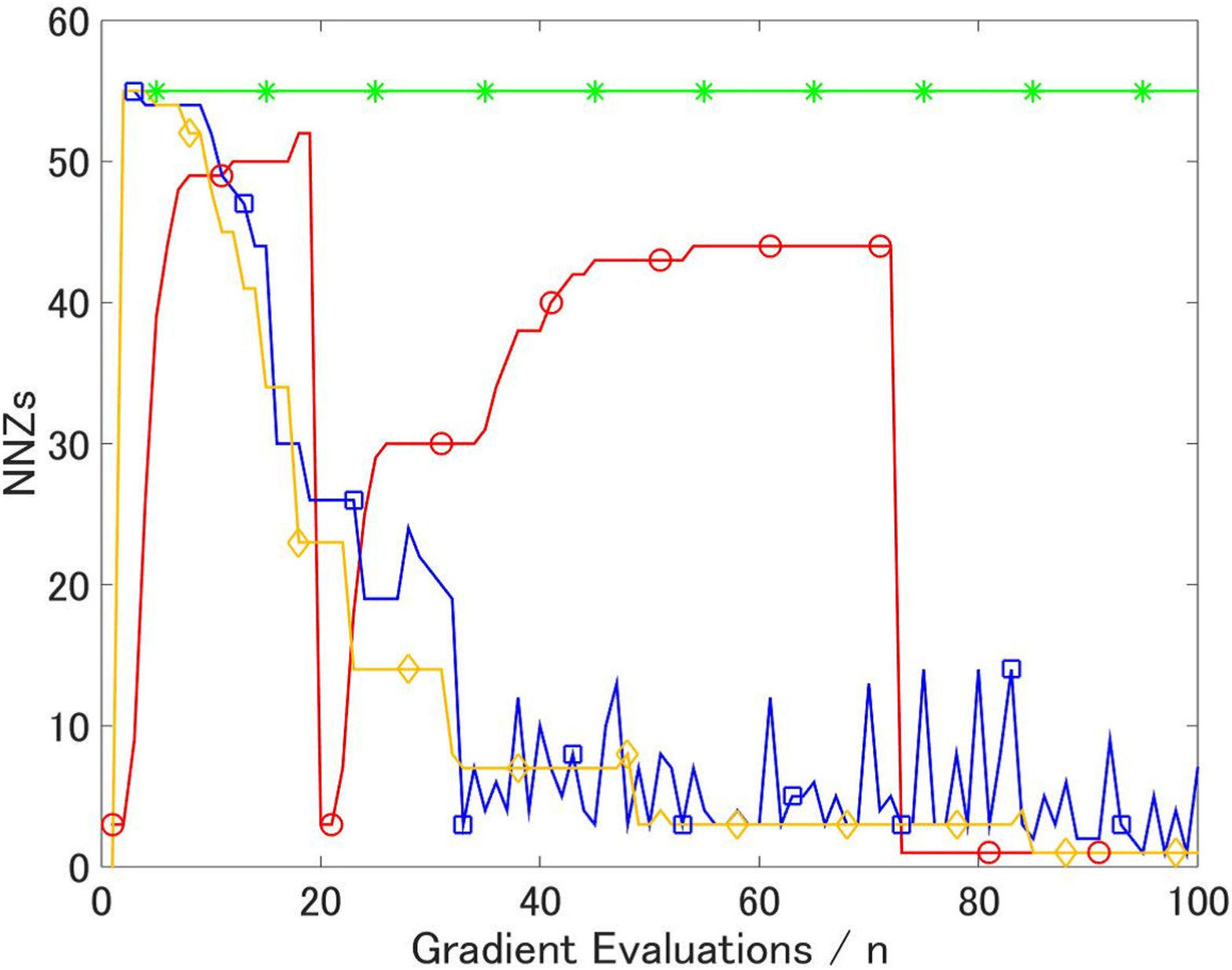}
\end{minipage}
}
\caption{Comparison of different methods on covertype data set}
\label{covertype}
\end{figure}
Figures \ref{covertype_sc} and \ref{covertype_nsc} show the comparison of SVRDA and SADA with the different methods described above on the covertype data set for different setups of $\lambda_1$ and $\lambda_2$ (the strongly convex case $\lambda_2=10^{-6} >0$ (Figure \ref{covertype_sc}) and the non-strongly convex case $\lambda_2=0$ (Figure \ref{covertype_nsc})). Objective Gap (left) means $P(x)-P(x_*)$ for the output solution $x$ and NNZs (right) means the number of nonzeros in the output solution. SVRDA-x and SADA-x output the solution generated by the gradient descent update $\widetilde{x}_s$, and SVRDA-v and SADA-v output the one generated by the dual averaging update $\widetilde{v}_s$. We do not report SVRDA-v and SADA-v for a non-strongly convex regularizer, because it has no theoretical convergence guarantee. For a strongly convex regularizer (top), UniVR, SVRDA-x, and SVRDA-v outperform other methods, as indicated by the theories (see Table \ref{table}). Observe that the objective gaps of SVRDA-x and SVRDA-v (respectively SADA-x and SADA-v) are very close though SVRDA-v (respectively SADA-v) has no theoretical guarantee for convergence of the objective gap. Note that SVRDA-v (respectively SADA-v) gives sparser solutions than SVRDA-x (respectively SADA-x) and the other methods. For a non-strongly convex regularizer (bottom), SVRDA-x and SADA-x converge more quickly than both UniVR and SAGA. The sparsity pattern of the output solutions of UniVR is unstable and that of SAGA is very poor, because UniVR and SAGA need to average the history of the solutions and the averaged solutions could be nonsparse. In contrast, SVRDA-x and SADA-x show a nice sparsity recovery performance. \par  
\begin{figure}[t]
\subfloat[$\lambda_1=10^{-4}$, $\lambda_2=10^{-4}$]
{\label{Reuters_sc}
\begin{minipage}{0.5\hsize}
\centering 
\includegraphics[bb=0 0 1500 1000, width=9.5cm]{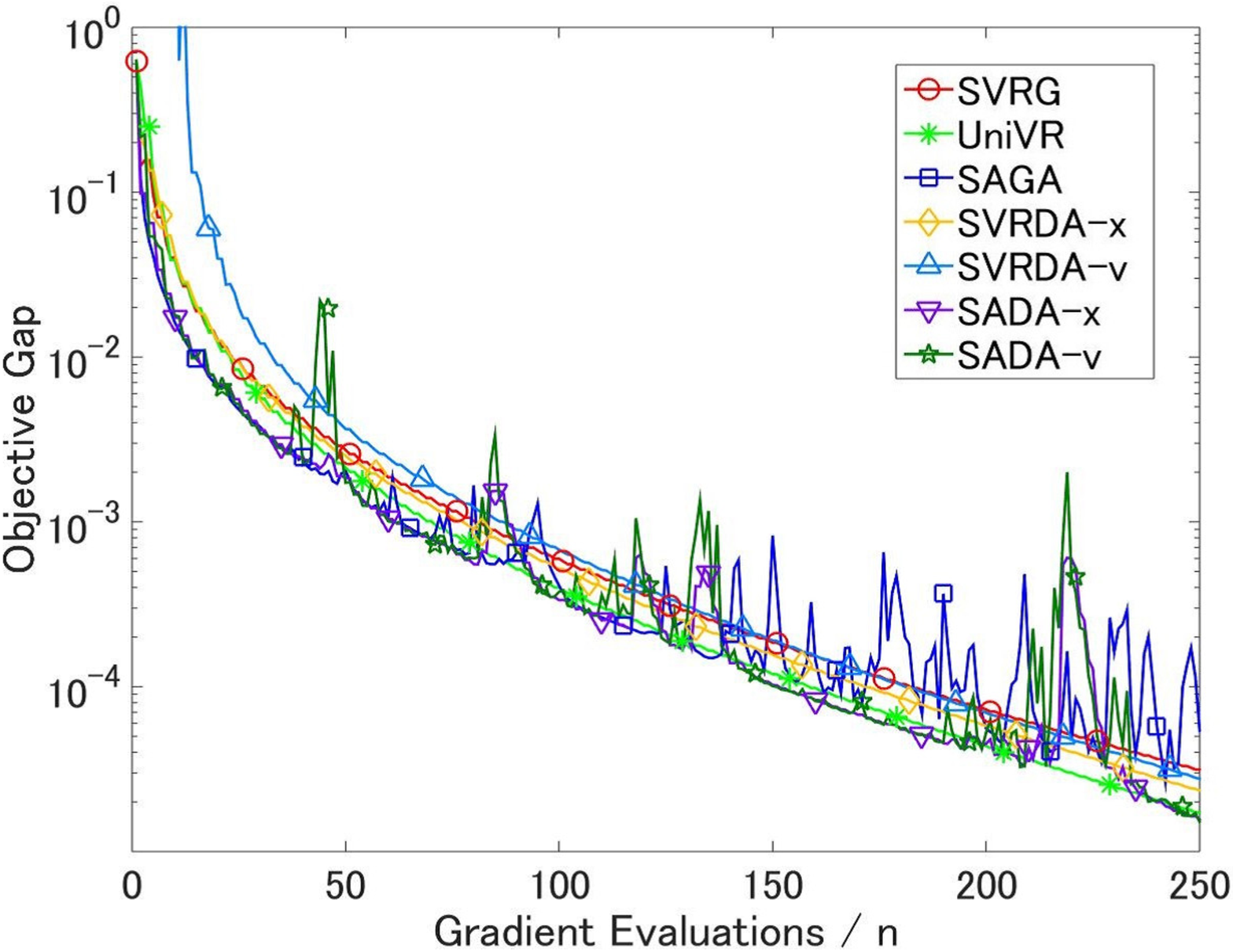}
\end{minipage}
\begin{minipage}{0.5\hsize}
\centering
\includegraphics[bb=0 0 1500 1000, width=9.5cm]{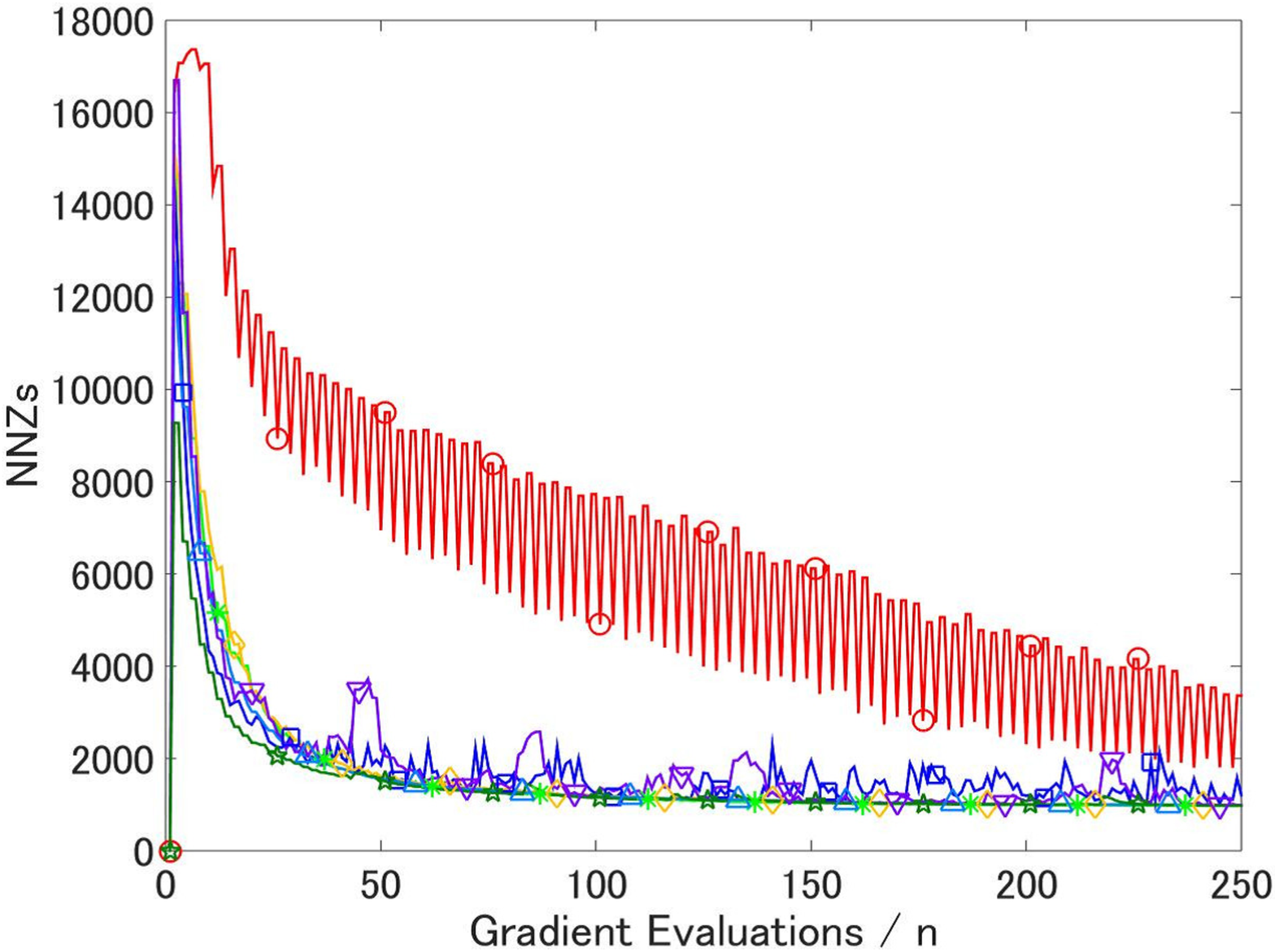}
\end{minipage}
}
\\
\subfloat[$\lambda_1=10^{-4}$, $\lambda_2=0$]{\label{Reuters_nsc}
\begin{minipage}{0.5\hsize}
\centering
\includegraphics[bb=0 0 1500 1000, width=9.5cm]{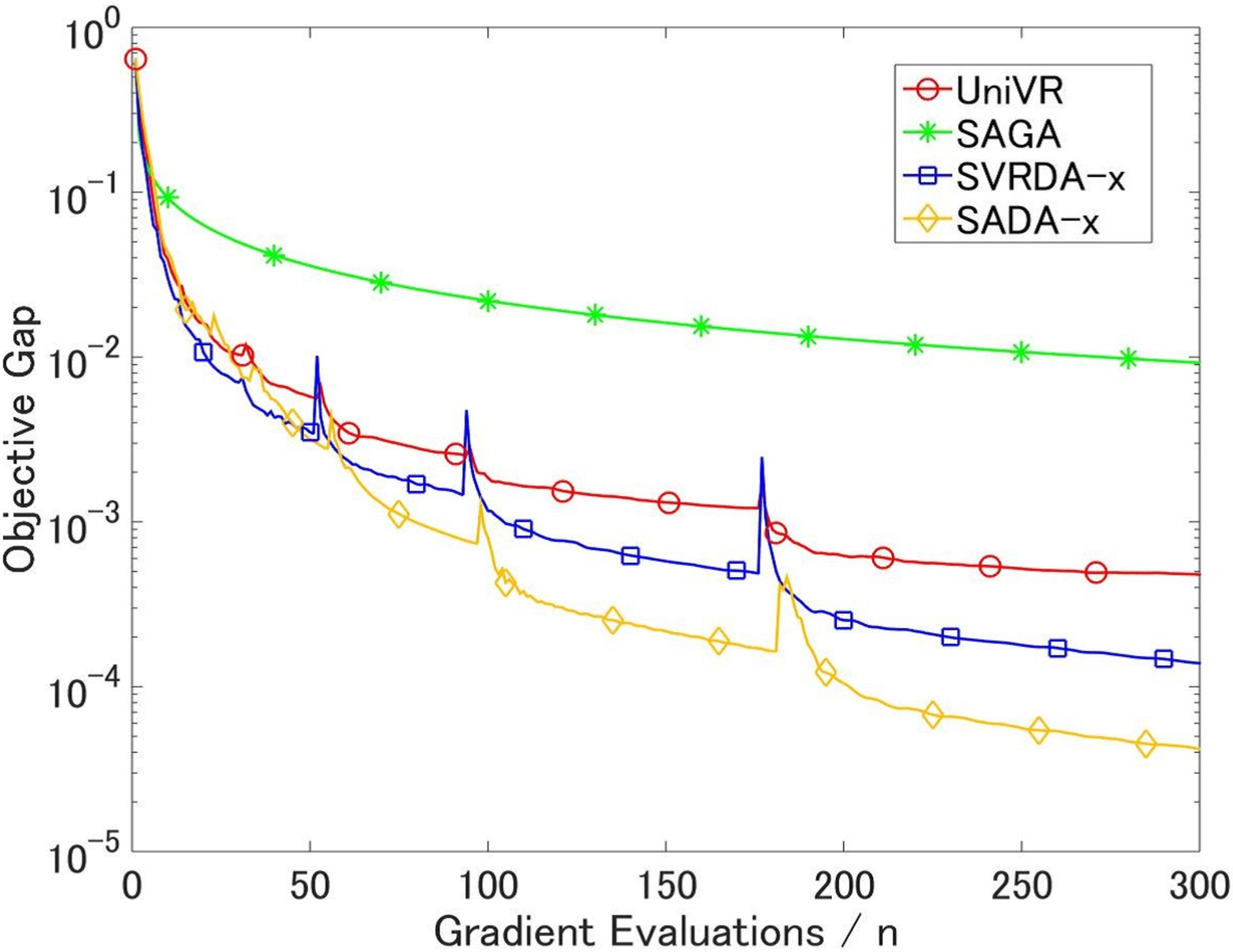}
\end{minipage}
\begin{minipage}{0.5\hsize}
\centering
\includegraphics[bb=0 0 1500 1000, width=9.5cm]{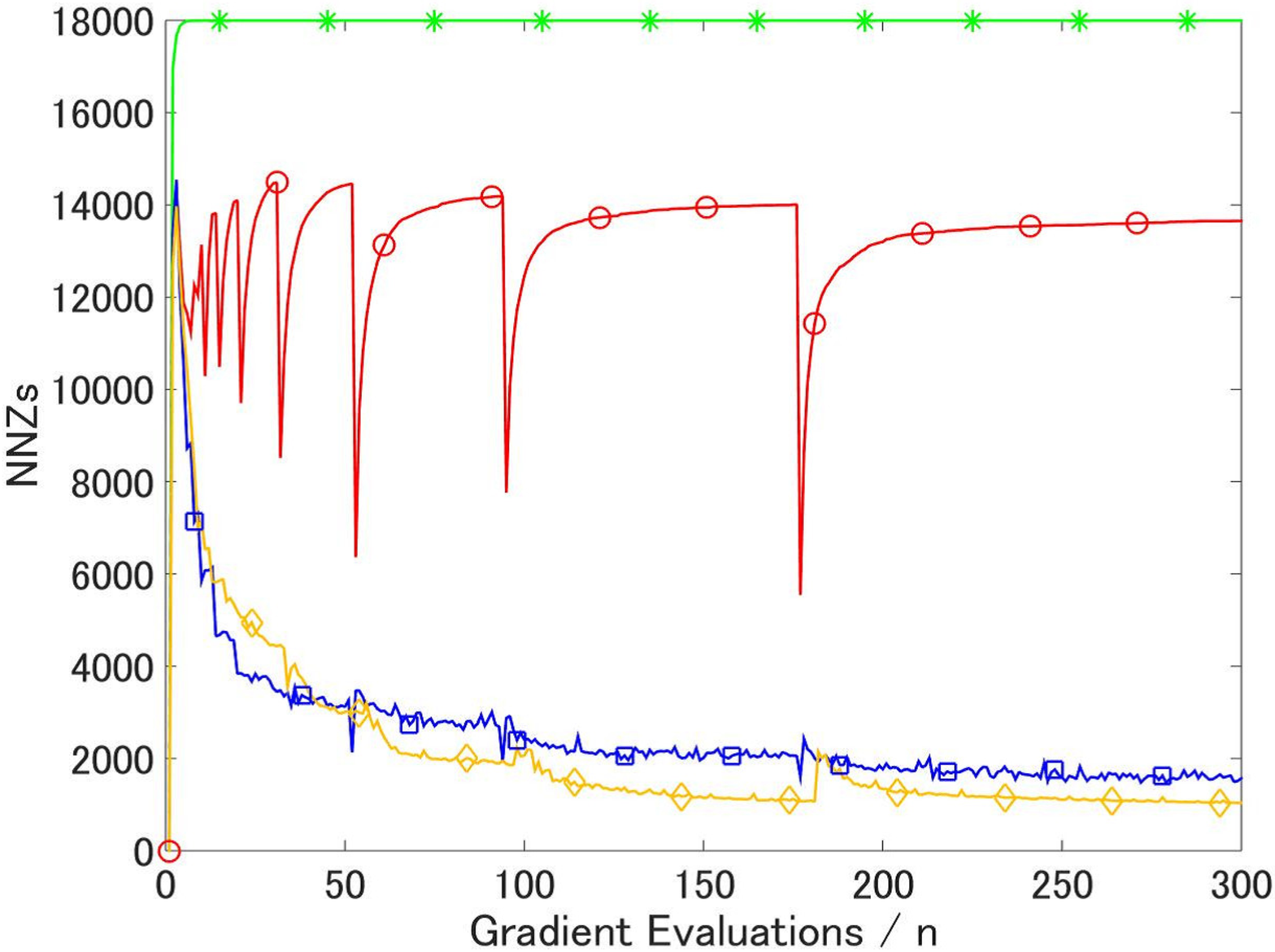}
\end{minipage}
}
\caption{Comparison of different methods on Reuters-21578 data set}
\label{Reuters}
\end{figure}
Figures \ref{Reuters_sc} and \ref{Reuters_nsc} show the comparison of different methods on the Reuters-21578 data set for different setups of $\lambda_1$ and $\lambda_2$ (the strongly convex case $\lambda_2=10^{-4} >0$ (Figure \ref{Reuters_sc}) and the non-strongly convex case $\lambda_2=0$ (Figure \ref{Reuters_nsc})). For a strongly convex regularizer, SVRG type algorithms (SVRG, UniVR, SVRDA-x, and SVRDA-v) show nice convergence behavior whereas SAGA type algorithms (SAGA, SADA-x, and SADA-v) show a slightly unstable behavior. Note that the sparsity pattern of the output solution of SVRG is poor. For a non-strongly convex regularizer, SVRDA-x and SADA-x converge more quickly but a bit more unstably than the other methods. Observe that, when a new stage starts, SVRDA-x and SADA-x lead to a sharp increase in the objective gap followed by a quick drop. This behavior can also be seen in the Multi-stage ORDA \citep{chen2012optimal}. We can see that the sparsity recovery performances of SVRDA-x and SADA-x are very nice whereas that of UniVR is unstable and poor and that of SAGA is quite poor. \par
\begin{figure}[t]
\subfloat[$\lambda_1=10^{-4}$, $\lambda_2=10^{-4}$]
{\label{sido0_sc}
\begin{minipage}{0.5\hsize}
\centering
\includegraphics[bb=0 0 1500 1000, width=9.5cm]{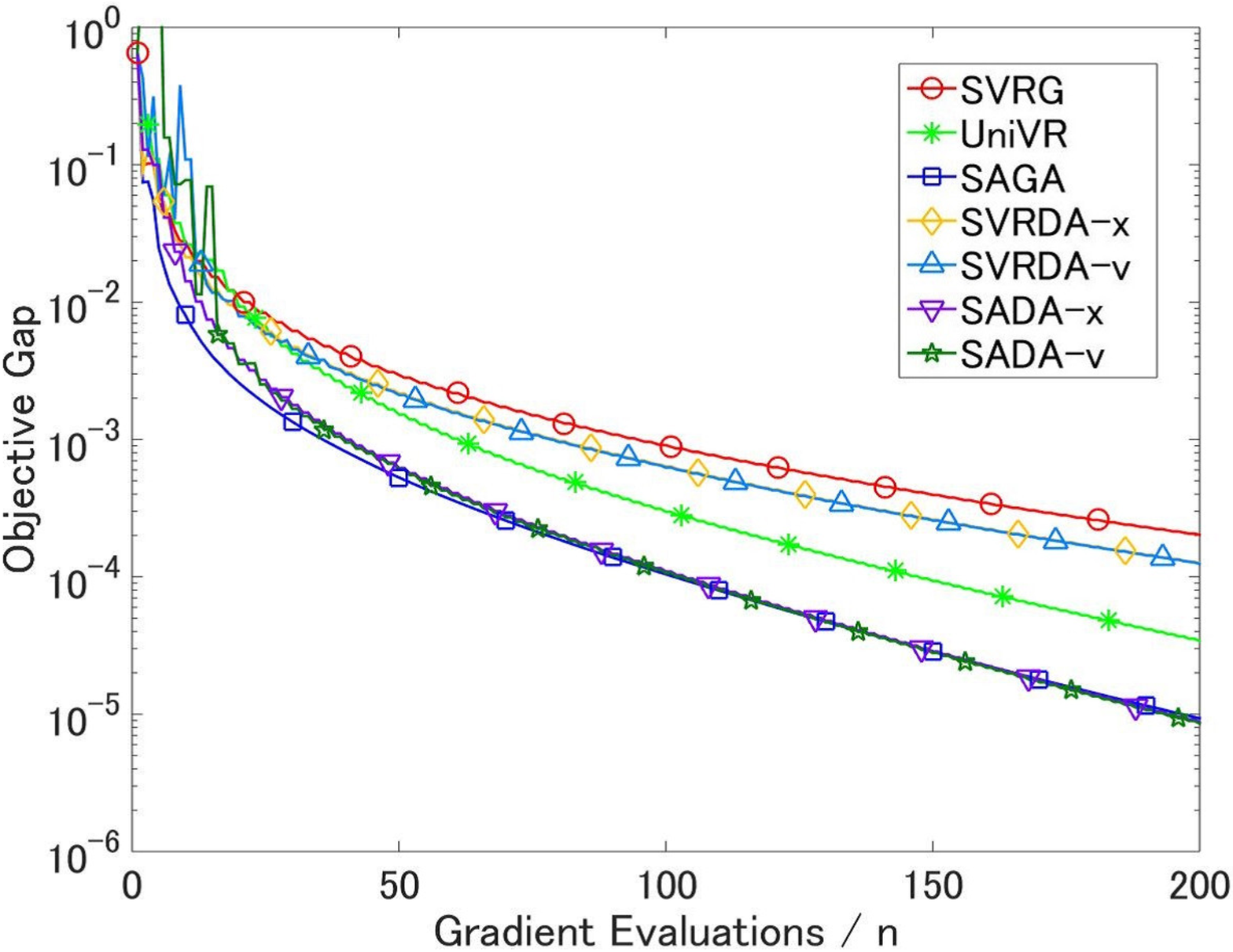}
\end{minipage}
\begin{minipage}{0.5\hsize}
\centering
\includegraphics[bb=0 0 1500 1000, width=9.5cm]{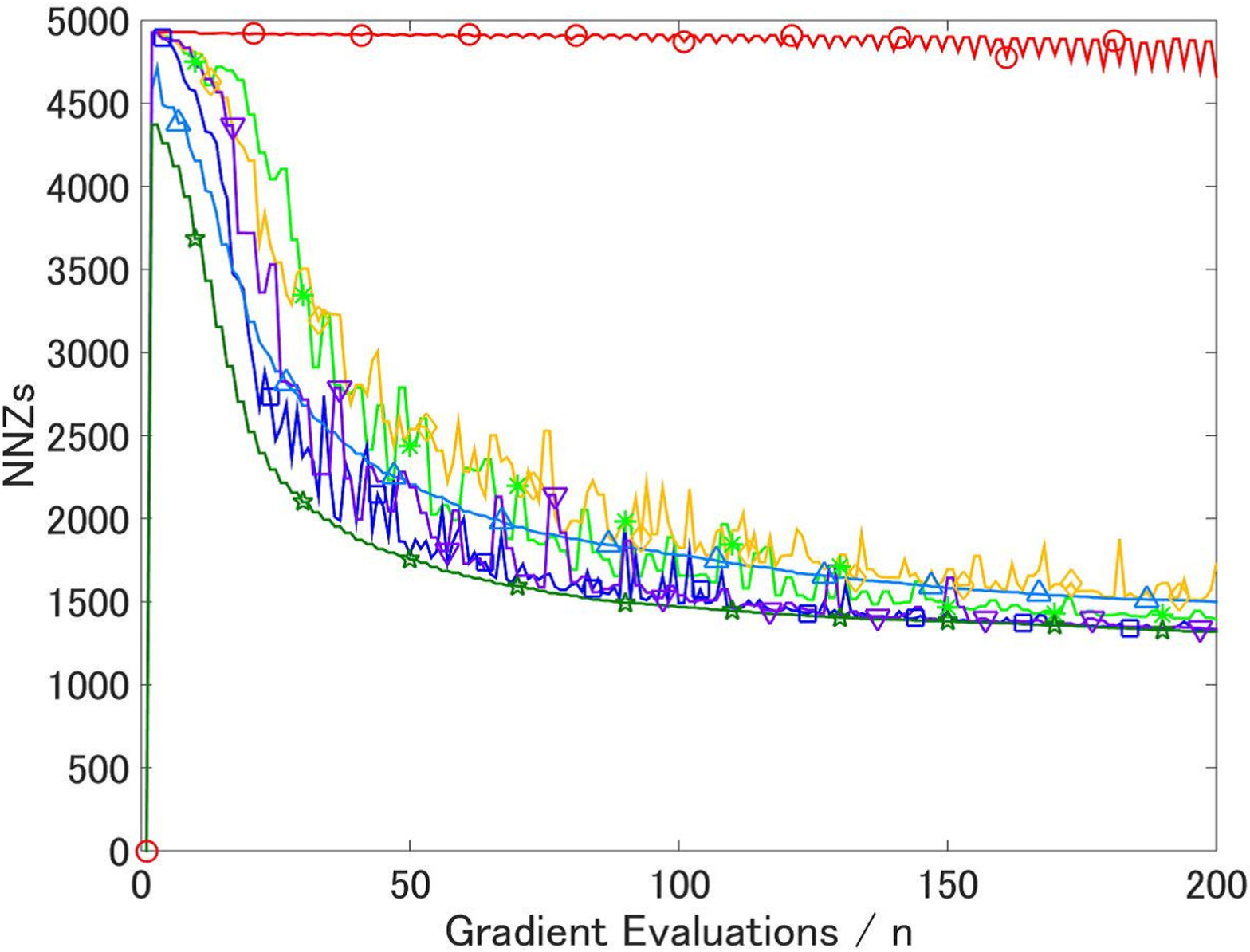}
\end{minipage}
}
\\
\subfloat[$\lambda_1=10^{-4}$, $\lambda_2=0$]{\label{sido0_nsc}
\begin{minipage}{0.5\hsize}
\centering
\includegraphics[bb=0 0 1500 1000, width=9.5cm]{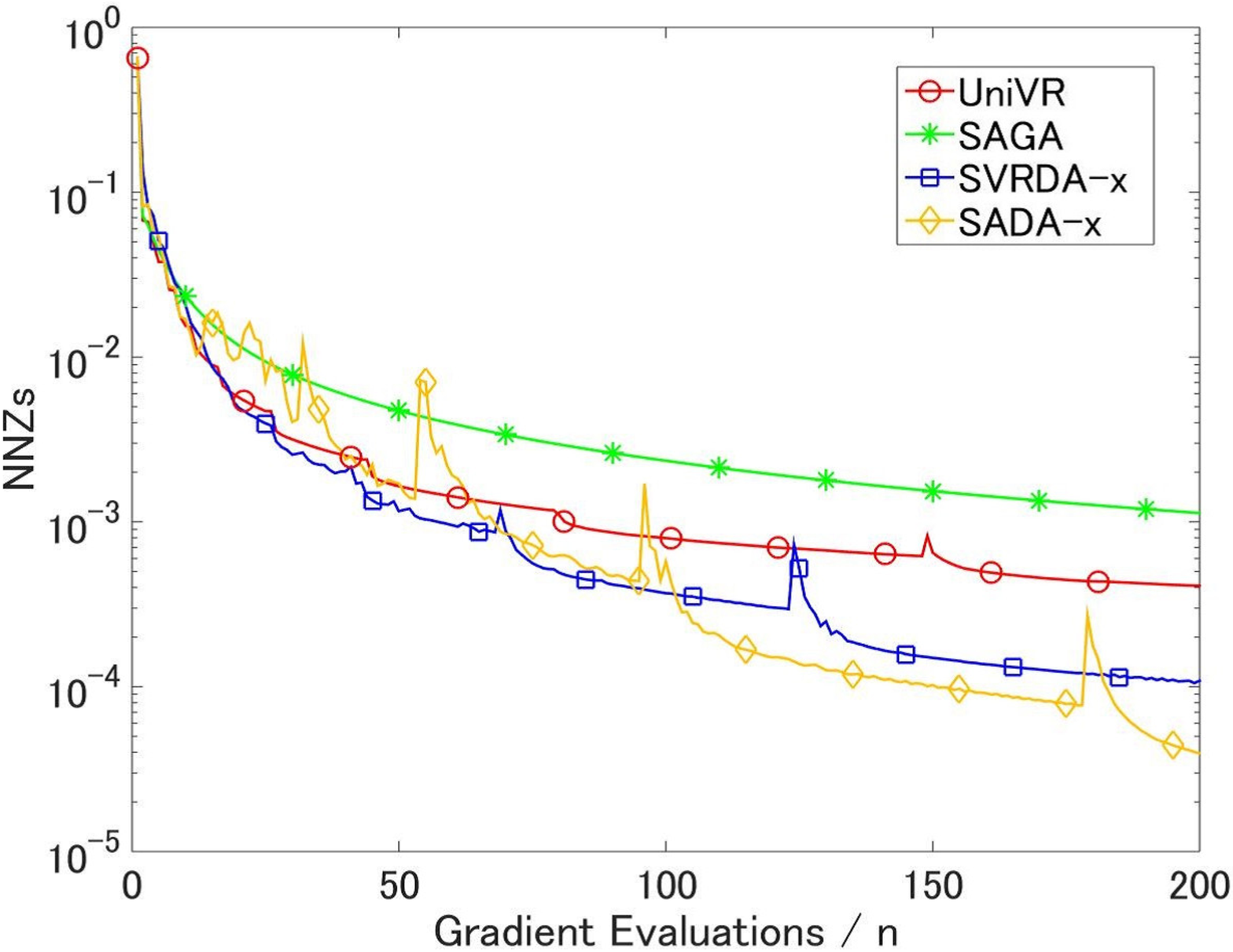}
\end{minipage}
\begin{minipage}{0.5\hsize}
\centering
\includegraphics[bb=0 0 1500 1000, width=9.5cm]{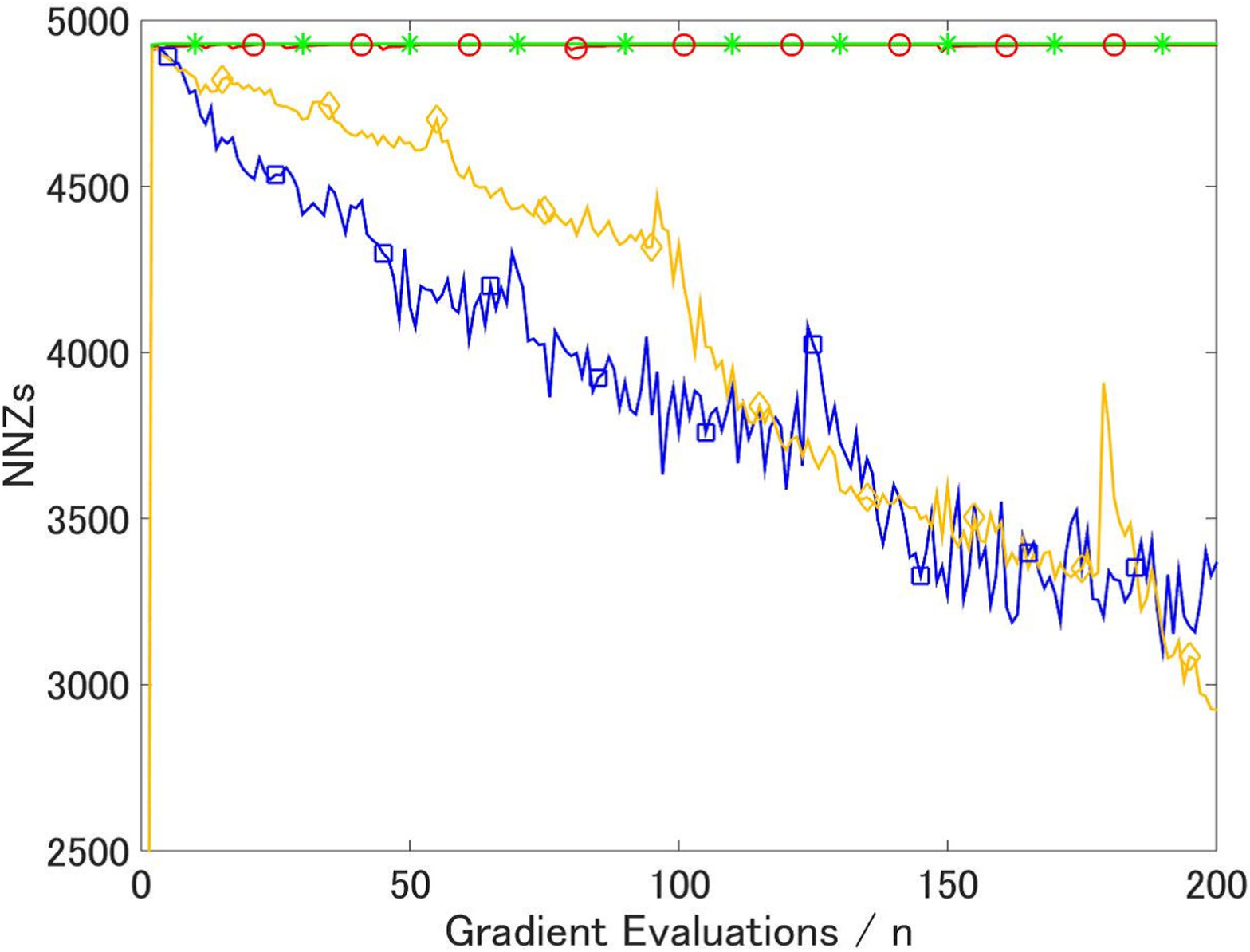}
\end{minipage}
}
\caption{Comparison of different methods on sido0 data set}
\label{sido0}
\end{figure}
Figures \ref{sido0_sc} and \ref{sido0_nsc} show the comparison of different methods on the sido0 data set for different setups of $\lambda_1$ and $\lambda_2$ (the strongly convex case $\lambda_2=10^{-4} >0$ (Figure \ref{sido0_sc}) and the non-strongly convex case $\lambda_2=0$ (Figure \ref{sido0_nsc})). For a strongly convex regularizer, the performances of SAGA, SADA-x, and SADA-v are among the best. Especially, SADA-v shows the best sparsity recovery performance. Note that the sparsity recovery performance of SVRG is very poor. We can see that the convergence of the NNZs of SVRDA-v (respectively SADA-v) is superior to SVRDA-x (respectively SADA-x). For a non-strongly convex regularizer, SVRDA-x and SADA-x outperform both UniVR and SAGA. Especially, SVRDA-x and SADA-x show nice sparsity recovery performances though the solutions of UniVR and SAGA are not sparse at all.  
\section{Conclusion and future work}
In this paper, we proposed two stochastic gradient methods for regularized empirical risk minimization problems: SVRDA and SADA. We have shown that SVRDA and SADA achieve $O\left(\left(n+\frac{\bar L}{\mu}\right)\mathrm{log}\frac{1}{\varepsilon}\right)$ and $O\left(\left(n+\frac{L_{\mathrm{max}}}{\mu}\right)\mathrm{log}\frac{1}{\varepsilon}\right)$complexity, respectively, for a strongly convex regularizer and $O\left(n\mathrm{log}\frac{1}{\varepsilon}+\frac{\bar L}{\varepsilon}\right)$ and $O\left(n\mathrm{log}\frac{1}{\varepsilon}+\frac{L_{\mathrm{max}}}{\varepsilon}\right)$ complexity, respectively, for a non-strongly convex regularizer. \par
In numerical experiments, our methods led to better sparsity recovery than the existing methods for sparsity-inducing regularizers and showed nice convergence behaviors, especially for non-strongly convex regularizers. \par 
An interesting future work is to extend our methods to the alternating directional multiplier method (ADMM) framework. In this paper, we assumed that the proximal mapping of $R$ can be efficiently computed. However, for structured regularization problems (for example, overlapped group lasso, graph lasso, etc.), this assumption is generally not satisfied and our methods cannot be directly applied. In contrast, ADMM can be applied to these problems without this assumption. Suzuki \citep{suzuki2013dual} has proposed regularized dual averaging-ADMM (RDA-ADMM), which is RDA \citep{xiao2009dual} for ADMM in an online setting. Furthermore, Suzuki \citep{suzuki2014stochastic} has proposed stochastic dual coordinate ascent-ADMM (SDCA-ADMM), which is SDCA \citep{shalev2013stochastic, shalev2013accelerated} for ADMM in regularized an empirical risk minimization setting,  and has shown that it converges exponentially for a strongly convex regularizer. Applying SVRDA to the ADMM framework and showing linear convergence for a strongly convex regularizer would be  promising future work. 

\section*{Acknowledgement}
This work was partially supported by MEXT Kakenhi (25730013, 25120012, and 26280009), JST-PRESTO and JST-CREST.
\bibliography{bibsvrda}
\newpage
\appendix
\title{Appendix}
\renewcommand\appendixname{}
\section{Proof of Theorem \ref{main thm1}} \label{appendix}
In this section, we give the proof of Theorem \ref{main thm1}. First we prove the following two easy lemmas.

\begin{lem}\label{avg lemma}
$$\bar g_t = \frac{1}{t}\sum_{\tau = 1}^{t} g_\tau \ (t \geq 1).$$
\end{lem}
\begin{pr}
For $t=1$, $\bar g_1 = g_1 = \frac{1}{1}\sum_{\tau = 1}^1 g_{\tau}$. \\
Assume that the claim holds for some $t \geq 1$. Then
\begin{align}
\bar g_{t+1} &= \left(1-\frac{1}{t+1}\right)\bar g_t + \frac{1}{t+1}g_{t+1} \ \text{(by the definition)} \notag \\
&= \left(1-\frac{1}{t+1}\right)\frac{1}{t}\sum_{\tau = 1}^{t} g_\tau + \frac{1}{t+1}g_{t+1} \ \text{(by the assumption of the induction)} \notag \\
&= \frac{1}{t+1}\sum_{\tau = 1}^{t+1} g_\tau. \notag
\end{align}
This finishes the proof for Lemma \ref{avg lemma}. \qed
\end{pr}

\begin{lem}\label{smooth lemma}
For every $x$, $u \in \mathbb{R}^d$, 
$$F(u) + \langle \nabla F(u), x-u \rangle + R(x)  \leq P(x) - \frac{1}{2\bar L}\frac{1}{n}\sum_{i=1}^{n}\frac{1}{nq_i}||\nabla f_i(x) - \nabla f_i(u)||^2.$$
\end{lem}
\begin{pr}
Since $f_i$ is $L_i$-smooth, we have (see \citep{nesterov2013introductory})
$$f_i(u) + \langle \nabla f_i(u), x-u \rangle \leq f_i (x) - \frac{1}{2L_i}||\nabla f_i(x) - \nabla f_i(u)||^2. $$
Summing this inequality from $i=1$ to $n$ and dividing it by $n$ results in 
$$F(u) + \langle \nabla F(u), x-u \rangle \leq F(x) - \frac{1}{2\bar L}\frac{1}{n}\sum_{i=1}^{n}\frac{1}{nq_i}||\nabla f_i(x) - \nabla f_i(u)||^2.$$
Adding $R(x)$ gives the desired result. \qed
\end{pr}

Next we prove the following main lemma.
\begin{lem}\label{main lemma1}
For the $s$th stage of SVRDA, 
\begin{align}
&\mathbb{E} [P({x}_{m_s})-P(x
_{*})] + \frac{\eta+m_s\mu}{2m_s}  \mathbb{E}||{v}_{m_s} - x_*||^2 \notag \\
\leq &\frac{1}{2 m_s} \sum_{t=1}^{m_s} \left[ \frac{t}{\eta t - \bar L}\mathbb{E} ||g_t- \nabla F(u_{t-1})||^2 -\frac{1}{\bar L}\mathbb{E} \left[\frac{1}{n}\sum_{i=1}^{n}\frac{1}{nq_i}||\nabla f_i(u_{t-1}) - \nabla f_i(x_*)||^2\right]\right] \notag \\
&+\frac{\eta}{2m_s}||{v}_{0} - x_*||^2, \notag
\end{align}
where the expectations are conditioned on all previous stages. 
\end{lem}
\begin{pr}
First note that $u_t = \left(1-\frac{1}{t+1}\right)x_t + \frac{1}{t+1}v_t$ for $t \geq 0$ by the definition of $u_0$. We define 
\begin{align}
\ell_t (x) &= F(u_{t-1}) + \langle \nabla F(u_{t-1}), x-u_{t-1} \rangle + R(x), \notag \\
\hat \ell_t (x) &= F(u_{t-1}) + \langle g_t, x-u_{t-1} \rangle + R(x). \notag
\end{align}
Observe that $\ell_t  \leq P$. For $t \geq 1$, by  Lemma \ref{avg lemma}, we have
\begin{align}
\sum_{\tau=1}^{t} \hat \ell_{\tau} (x) =& \sum_{\tau = 1}^{t} F(u_{t-1}) + \sum_{\tau = 1}^{t}\langle g_{\tau}, x-u_{\tau-1} \rangle + \sum_{\tau = 1}^{t} R(x) \notag \\
=& \langle t \bar g_t, x \rangle + t R(x) + \sum_{\tau = 1}^{t}F(u_{t-1}) - \sum_{\tau = 1}^{t}\langle g_{\tau}, u_{\tau-1} \rangle \notag
\end{align}
and thus we have $v_t = \underset{x \in \mathbb{R}^d} { \mathrm{argmin}\ }\left\{\sum_{\tau=1}^{t} \hat \ell_{\tau} (x) + \frac{\eta}{2}||x-v_0||^2 \right\} $. \\Also note that $x_t = \underset{x \in \mathbb{R}^d} { \mathrm{argmin}\ }\left\{ \hat \ell_{\tau} (x) + \frac{\eta t}{2}||x-u_{t-1}||^2 \right\} $. Since $F$ is $\bar L$-smooth, we have (see \citep{nesterov2013introductory})
$$F(x_t) \leq F(u_{t-1}) + \langle \nabla F(u_{t-1}), x_t - u_{t-1} \rangle + \frac{\bar L}{2} ||x_t-u_{t-1}||^2,$$and thus 
\begin{align}
P(x_t) \leq& \ell_t(x_t) + \frac{\bar L}{2}||x_t-u_{t-1}||^2 \notag \\
=& \hat \ell_t(x_t) + \frac{\eta t}{2}||x_t-u_{t-1}||^2 - \frac{\eta t -\bar L}{2}||x_t-u_{t-1}||^2 -\langle g_t-\nabla F(u_{t-1}), x_t- u_{t-1} \rangle. \notag
\end{align}
Since $x_t$ is the minimizer of $\hat \ell_t (x) + \frac{\eta t}{2}||x-u_{t-1}||^2$, we have
\begin{align}
\hat \ell_t (x_t) +& \frac{\eta t}{2}||x_t-u_{t-1}||^2 \notag \\ \leq&\hat \ell_t \left(\left(1-\frac{1}{t}\right)x_{t-1} + \frac{1}{t}v_t\right) + \frac{\eta t}{2}\left|\left|\left(1-\frac{1}{t}\right)x_{t-1} + \frac{1}{t}v_t-u_{t-1}\right|\right|^2 \notag
\end{align}
and hence 
\begin{align}
P(x_t) \leq& 
\hat \ell_t \left(\left(1-\frac{1}{t}\right)x_{t-1} + \frac{1}{t}v_t\right) + \frac{\eta t}{2}\left|\left|\left(1-\frac{1}{t}\right)x_{t-1} + \frac{1}{t}v_t-u_{t-1}\right|\right|^2 \notag \\
&- \frac{\eta t - \bar L}{2}||x_t-u_{t-1}||^2- \langle g_t-\nabla F(u_{t-1}), x_t- u_{t-1} \rangle. \notag
\end{align}
Using the convexity of $\hat \ell_t$ and the  facts that $\left(1-\frac{1}{t}\right)x_{t-1} + \frac{1}{t}v_t-u_{t-1} = \frac{1}{t}(v_t-v_{t-1})$ and
\begin{align} - \frac{\eta t - \bar L}{2}||x_t-u_{t-1}||^2-& \langle g_t-\nabla F(u_{t-1}), x_t- u_{t-1} \rangle  \leq \frac{1}{2(\eta t - \bar L)}||g_t- \nabla F(u_{t-1})||^2,\notag 
\end{align} we get
\begin{align}
P(x_t) \leq& \left(1-\frac{1}{t}\right)\hat \ell_t (x_{t-1}) + \frac{1}{t}\hat \ell_t (v_t) +\frac{\eta}{2t}||v_t -v_{t-1}||^2 
+ \frac{1}{2(\eta t - \bar L)}||g_t- \nabla F(u_{t-1})||^2 \notag \\
=& \left(1-\frac{1}{t}\right) \ell_t (x_{t-1}) + \frac{1}{t}\hat \ell_t (v_t) +\frac{\eta}{2t}||v_t -v_{t-1}||^2 \notag \\
&+ \frac{1}{2(\eta t - \bar L)}||g_t- \nabla F(u_{t-1})||^2 
+ \left(1-\frac{1}{t}\right)\langle g_t-\nabla F(u_{t-1}), x_{t-1}- u_{t-1} \rangle \notag \\
\leq&  \left(1-\frac{1}{t}\right)P(x_{t-1})  + \frac{1}{t}\hat \ell_t (v_t) +\frac{\eta}{2t}||v_t -v_{t-1}||^2 \notag \\ &+   \frac{1}{2(\eta t - \bar L)}||g_t- \nabla F(u_{t-1})||^2
+ \left(1-\frac{1}{t}\right)\langle g_t-\nabla F(u_{t-1}), x_{t-1}- u_{t-1} \rangle. \notag
\end{align}
Multiplying both sides of the above inequality by $t$, we have
\begin{align}
t P(x_t) \leq& (t-1)P(x_{t-1}) + \hat \ell_t (v_t) + \frac{\eta}{2}||v_t -v_{t-1}||^2 \notag \\ &+ \frac{t}{2(\eta t - \bar L)}||g_t- \nabla F(u_{t-1})||^2 
+ (t-1)\langle g_t-\nabla F(u_{t-1}), x_{t-1}- u_{t-1} \rangle. \notag
\end{align}
By the fact that $\sum_{\tau=1}^{t-1} \hat \ell_{\tau} (x) + \frac{\eta}{2}||x-v_0||^2$ is $\eta$-strongly convex and $v_{t-1}$ is the minimizer of $\sum_{\tau=1}^{t-1} \hat \ell_{\tau} (x) +\frac{\eta}{2}||x-v_0||^2$ for $t \geq 2$, we have
$$\sum_{\tau=1}^{t-1} \hat \ell_{\tau} (v_{t-1}) +\frac{\eta}{2}||v_{t-1}-v_0||^2 +  \frac{\eta}{2}||v_t -v_{t-1}||^2 \leq\sum_{\tau=1}^{t-1} \hat \ell_{\tau} (v_t) + \frac{\eta}{2}||v_t-v_0||^2$$for $t \geq 1$ (and, for $t=1$, we define $\sum_{\tau=1}^0 = 0$). Using this inequality, we obtain 
\begin{align}
&t P(x_t) - \sum_{\tau=1}^{t} \hat \ell_{\tau} (v_t) - \frac{\eta}{2}||v_t-v_0||^2 \notag \\\leq& (t-1)P(x_{t-1}) -\sum_{\tau=1}^{t-1} \hat \ell_{\tau} (v_{t-1}) - \frac{\eta}{2}||v_{t-1}-v_0||^2
+ \frac{t}{2(\eta t - \bar L)}||g_t- \nabla F(u_{t-1})||^2 \notag \\ &+ (t-1)\langle g_t-\nabla F(u_{t-1}), x_{t-1}- u_{t-1} \rangle. \notag
\end{align}
Summing the above inequality from $t=1$ to $m_s$ results in
\begin{align}
m_s P(x_{m_s}) - &\sum_{t=1}^{m_s} \hat \ell_t (v_{m_s}) - \frac{\eta}{2}||v_{m_s}-v_0||^2 \notag \\
\leq& \sum_{t=1}^{m_s}\frac{t}{2(\eta t - \bar L)}||g_t- \nabla F(u_{t-1})||^2 \notag \\
&+ \sum_{t=1}^{m_s}(t-1)\langle g_t-\nabla F(u_{t-1}), x_{t-1}- u_{t-1} \rangle. \notag
\end{align}
Using $\eta+m_s\mu$-strongly convexity of  the function $\sum_{t=1}^{m_s} \hat \ell_t (x) + \frac{\eta}{2}||x-v_0||^2$ and the optimality of $v_{m_s}$, we have
$$\sum_{t=1}^{m_s} \hat \ell_t (v_{m_s}) + \frac{\eta}{2}||v_{m_s}-v_0||^2 \leq \sum_{t=1}^{m_s} \hat \ell_t (x_*) + \frac{\eta}{2}||v_0-x_*||^2 - \frac{\eta +m_s \mu}{2}||v_{m_s}-x_*||^2$$
and hence
\begin{align}
&m_s P(x_{m_s}) \notag \\ \leq&  \sum_{t=1}^{m_s} \hat \ell_t (x_*) + \frac{\eta}{2}||v_0-x_*||^2 - \frac{\eta +m_s \mu}{2}||v_{m_s}-x_*||^2 \notag \\
&+\sum_{t=1}^{m_s}\frac{t}{2(\eta t - \bar L)}||g_t- \nabla F(u_{t-1})||^2 + \sum_{t=1}^{m_s}(t-1)\langle g_t-\nabla F(u_{t-1}), x_{t-1}- u_{t-1} \rangle \notag \\
=&\sum_{t=1}^{m_s} \ell_t (x_*) + \frac{\eta}{2}||v_0-x_*||^2 - \frac{\eta +m_s \mu}{2}||v_{m_s}-x_*||^2 \notag \\
&+\sum_{t=1}^{m_s}\frac{t}{2(\eta t - \bar L)}||g_t- \nabla F(u_{t-1})||^2+ \sum_{t=1}^{m_s}(t-1)\langle g_t-\nabla F(u_{t-1}), x_{t-1}- u_{t-1} \rangle \notag \\
&+\sum_{t=1}^{m_s}\langle g_t-\nabla F(u_{t-1}), x_*- u_{t-1} \rangle. \notag 
\end{align}
By Lemma \ref{smooth lemma} with $x = x_*$ and $u = u_{t-1}$, we have
$$\ell_t(x_*)\leq P(x_*) - \frac{1}{2\bar L}\frac{1}{n}\sum_{i=1}^{n}\frac{1}{nq_i}||\nabla f_i(x_*) - \nabla f_i(u_{t-1})||^2.$$
Applying this inequality to the above inequality yields
\begin{align}
&m_s P(x_{m_s})\notag \\ \leq&m_s P(x_*) +\frac{\eta}{2}||v_0-x_*||^2 - \frac{\eta +m_s \mu}{2}||v_{m_s}-x_*||^2 \notag \\
&+ \frac{1}{2}\sum_{t=1}^{m_s} \left[ \frac{t}{\eta t - \bar L}||g_t- \nabla F(u_{t-1})||^2 -\frac{1}{\bar L}\frac{1}{n}\sum_{i=1}^{n}\frac{1}{nq_i}||\nabla f_i(x_*) - \nabla f_i(u_{t-1})||^2\right] \notag \\
&+ \sum_{t=1}^{m_s}(t-1)\langle g_t-\nabla F(u_{t-1}), x_{t-1}- u_{t-1} \rangle 
+\sum_{t=1}^{m_s}\langle g_t-\nabla F(u_{t-1}), x_*- u_{t-1} \rangle. \notag 
\end{align}
Dividing this inequality by $m_s$ results in
\begin{align}
&P(x_{m_s})\notag \\  \leq& P(x_*)+ \frac{\eta}{2m_s}||v_0-x_*||^2 - \frac{\eta +m_s \mu}{2m_s }||v_{m_s}-x_*||^2 \notag \\
&+\frac{1}{2 m_s} \sum_{t=1}^{m_s} \left[ \frac{t}{\eta t - \bar L}||g_t- \nabla F(u_{t-1})||^2 -\frac{1}{\bar L}\frac{1}{n}\sum_{i=1}^{n}\frac{1}{nq_i}||\nabla f_i(x_*) - \nabla f_i(u_{t-1})||^2\right] \notag \\
&+ \frac{1}{m_s}\sum_{t=1}^{m_s}(t-1)\langle g_t-\nabla F(u_{t-1}), x_{t-1}- u_{t-1} \rangle 
+\frac{1}{m_s}\sum_{t=1}^{m_s}\langle g_t-\nabla F(u_{t-1}), x_*- u_{t-1} \rangle. \notag 
\end{align}
Taking the expectation on both sides yields
\begin{align}
&\mathbb{E} [P({x}_{m_s})-P(x
_{*})] + \frac{\eta+m_s\mu}{2m_s}  \mathbb{E}||{v}_{m_s} - x_*||^2 \notag \\
\leq &\frac{1}{2 m_s} \sum_{t=1}^{m_s} \left[ \frac{t}{\eta t - \bar L}\mathbb{E} ||g_t- \nabla F(u_{t-1})||^2 -\frac{1}{\bar L}\mathbb{E} \left[\frac{1}{n}\sum_{i=1}^{n}\frac{1}{nq_i}||\nabla f_i(x_*) - \nabla f_i(u_{t-1})||^2\right]\right] \notag \\
&+\frac{\eta}{2m_s}||{v}_{0} - x_*||^2. \notag
\end{align}
Here we used the fact that $\mathbb{E}[ g_t-\nabla F(u_{t-1})] = 0$ for $t=1, \ldots , m_s$. This finishes the proof of Lemma \ref{main lemma1}. 
\qed
\end{pr}
Now we need the following lemma. 
\begin{lem}\label{object bound lemma}
For every $x \in \mathbb{R}^d$, 
$$\frac{1}{n}\sum_{i=1}^{n}\frac{1}{nq_i}||\nabla f_i(x) - \nabla f_i(x_*)||^2 \leq 2 \bar L (P(x) - P(x_*)-\frac{\mu}{2}||x-x_*||^2).$$
\end{lem}
\begin{pr}
From the argument of the proof of Lemma \ref{smooth lemma}, we have
$$\frac{1}{n}\sum_{i=1}^{n}\frac{1}{nq_i}||\nabla f_i(x) - \nabla f_i(x_*)||^2 \leq 2\bar L (F(x)-\langle \nabla F(x_*), x-x_* \rangle - F(x_*)).$$
By the optimality  of $x_*$, there exists $\xi_* \in \partial R(x_*)$ such that $\nabla F(x_*) + \xi_*$. Then using $\mu$-strong convexity of $R$, we get
$$-\langle \nabla F(x_*), x-x_* \rangle = \langle \xi_*, x-x_* \rangle \leq R(x)-R(x_*)-\frac{\mu}{2}||x-x_*||^2$$
and hence 
$$\frac{1}{n}\sum_{i=1}^{n}\frac{1}{nq_i}||\nabla f_i(x) - \nabla f_i(x_*)||^2 \leq 2 \bar L (P(x) - P(x_*)-\frac{\mu}{2}||x-x_*||^2).$$\qed
\end{pr}
\begin{pr main thm1}
We bound the term $\mathbb{E}||g_t - \nabla F(u_{t-1})||^2$: 
\begin{align}
&\mathbb{E}||g_t - \nabla F(u_{t-1})||^2 \notag \\
=& \mathbb{E}\left[ \mathbb{E}_{i_t}\left[||(\nabla f_{i_t}(u_{t-1})-\nabla f_{i_t} (x_0))/n q_{i_t}+\nabla F(x_0)-\nabla F(u_{t-1})||^2|i_1, \ldots , i_{t-1} \right]\right] \notag \\
\leq&\mathbb{E}\left[ \mathbb{E}_{i_t}\left[||(\nabla f_{i_t}(u_{t-1})-\nabla f_{i_t} (x_0))/n q_{i_t}||^2|i_1, \ldots , i_{t-1} \right]\right] \notag \\
=& \mathbb{E}\left[\frac{1}{n}\sum_{i=1}^n \frac{1}{n q_{i}}||\nabla f_{i}(u_{t-1})-\nabla f_{i} (x_0)||^2\right] \notag \\
\leq& 3\mathbb{E}\left[\frac{1}{n}\sum_{i=1}^n \frac{1}{n q_{i}}||\nabla f_{i}(u_{t-1})-\nabla f_{i} (x_*)||^2\right] + \frac{3}{2}\left[\frac{1}{n}\sum_{i=1}^n \frac{1}{n q_{i}}||\nabla f_{i}(x_0)-\nabla f_{i} (x_*)||^2\right]. \notag 
\end{align}
Combining this inequality with Lemme \ref{object bound lemma}, we get
\begin{align}
\mathbb{E}||g_t - \nabla F(u_{t-1})||^2 \leq&  3\mathbb{E}\left[\frac{1}{n}\sum_{i=1}^n \frac{1}{n q_{i}}||\nabla f_{i}(u_{t-1})-\nabla f_{i} (x_*)||^2\right] \notag \\
&+ 3\bar L (P(x_0) - P(x_*)-\frac{\mu}{2}||x_0-x_*||^2). \notag
\end{align}
Since $\eta = \frac{1}{4 \bar L}$, using the inequality 
\begin{equation}\frac{t}{\eta t - \bar L} \leq \frac{1}{3\bar L} \ \ \ \ (\forall t \geq 1), \label{roughbound}\end{equation}
by Lemma \ref{main lemma1} we obtain
\begin{align}
&\mathbb{E} [P({x}_{m_s})-P(x_{*})] + \frac{\eta+m_s\mu}{2m_s}  \mathbb{E}||{v}_{m_s} - x_*||^2\notag \\ \leq& \frac{1}{2}(P(x_0) - P(x_*)-\frac{\mu}{2}||x_0-x_*||^2)
+\frac{\eta}{2m_s}||{v}_{0} - x_*||^2. \notag 
\end{align}
Since $x_{m_s}=\widetilde{x}_s$, $v_{m_s}=\widetilde{v}_s$, $x_0=\widetilde{x}_{s-1}$, and  $v_0=(1-\alpha)\widetilde{v}_{s-1}+\alpha\widetilde{x}_{s-1}$, we have
\begin{align}
&\mathbb{E} [P(\widetilde{x}_{s})-P(x_{*})] + \frac{\eta+m_s\mu}{2m_s} \mathbb{E}||\widetilde{v}_{s} - x_*||^2\notag \\ 
\leq& \frac{1}{2} (P(\widetilde{x}_{s-1})-P(x_{*})) + \left(\frac{\alpha\eta }{2m_s}-\frac{\mu}{4}\right)||\widetilde{x}_{s-1} - x_*||^2+\frac{(1-\alpha)\eta }{2m_s}||\widetilde{v}_{s-1} - x_*||^2.\notag 
\end{align}
Finally, taking expectations with respect to all previous stages gives the desired result. \qed 
\end{pr main thm1}
\section{Proof of Theorem \ref{main thm2}}\label{appendix2}
In this section, we give the proof of Theorem \ref{main thm2}.
\begin{lem}
\label{main lemma2}
For the $s$th stage of SADA, 
\begin{align}
&\mathbb{E} [P({x}_{m_s})-P(x_{*})] + \frac{\eta+m_s\mu}{2m_s} \mathbb{E}||{v}_{m_s} - x_*||^2 \notag \\
\leq& \frac{1}{2 m_s} \sum_{t=1}^{m_s} \left[ \frac{t}{\eta t-L_{\mathrm{max}}}\mathbb{E}||g_t - \nabla F(u_{t-1})||^2 -  \frac{1}{L_{\mathrm{max}}} \mathbb{E}\left[\frac{1}{n}\sum_{i=1}^{n}|| \nabla f_i(u_{t-1})-\nabla f_i(x_*) ||^2\right] \right] \notag \\  &+ \frac{\eta}{2m_s}||{v}_{0} - x_*||^2. \notag
\end{align}
\end{lem}
The proof of Lemma \ref{main lemma2} is identical to the proof of Lemma \ref{main lemma1} and we omit it. 
\begin{pr main thm2}
First we bound the term $\mathbb{E}||g_t - \nabla f(u_{t-1})||^2$: 
\begin{align}
\mathbb{E}||g_t - \nabla F(u_{t-1})||^2 =& \mathbb{E}\left[ \mathbb{E}_{i_t}\left[||\nabla f_{i_t}(u_{t-1})-\nabla f_{i_t} (\phi_{i_t}^{t-1})+\frac{1}{n}\sum_{i=1}^n \nabla f_i(\phi_i^{t-1})-\nabla F(u_{t-1})||^2|i_1, \ldots , i_{t-1} \right]\right] \notag \\
\leq&\mathbb{E}\left[ \mathbb{E}_{i_t}\left[||\nabla f_{i_t}(u_{t-1})-\nabla f_{i_t} (\phi_{i_t}^{t-1})||^2|i_1, \ldots , i_{t-1} \right]\right] \notag \\
=& \mathbb{E}\left[\frac{1}{n}\sum_{i=1}^n ||(\nabla f_{i}(u_{t-1})-\nabla f_{i} (\phi_{i}^{t-1}))||^2\right] \notag \\
\leq& 2 \mathbb{E}\left[\frac{1}{n}\sum_{i=1}^n ||(\nabla f_{i}(u_{t-1})-\nabla f_{i} (x_*))||^2\right] + 2\mathbb{E}\left[\frac{1}{n}\sum_{i=1}^n ||(\nabla f_{i}(\phi_{i}^{t-1})-\nabla f_{i} (x_*))||^2\right]. \notag 
\end{align}
Next we bound the term $\mathbb{E}\left[\frac{1}{n}\sum_{i=1}^n ||(\nabla f_{i}(\phi_{i}^{t-1})-\nabla f_{i} (x_*))||^2\right]$ for $t\geq 1$:
\begin{align}
&\mathbb{E}\left[\frac{1}{n}\sum_{i=1}^n ||(\nabla f_{i}(\phi_{i}^{t-1})-\nabla f_{i} (x_*))||^2\right] \notag \\
=& \mathbb{E}\left[ \mathbb{E}_{i_{t-1}}\left[\frac{1}{n}\sum_{i=1}^n ||(\nabla f_{i}(\phi_{i}^{t-1})-\nabla f_{i} (x_*))||^2|i_1, \ldots , i_{t-2}\right]\right]\notag \\
=& \mathbb{E}\left[\frac{1}{n} \sum_{i=1}^n\mathbb{E}_{i_{t-1}}\left[ ||(\nabla f_{i}(\phi_{i}^{t-1})-\nabla f_{i} (x_*))||^2|i_1, \ldots , i_{t-2}\right]\right] \notag \\
=&\mathbb{E}\left[\frac{1}{n} \sum_{i=1}^n\left[ \frac{1}{n} ||(\nabla f_{i}(u_{t-2})-\nabla f_{i} (x_*))||^2 + \left(1-\frac{1}{n}\right) ||(\nabla f_{i}(\phi_{i}^{t-2})-\nabla f_{i} (x_*))||^2\right] \right] \notag \\
=& \frac{1}{n}\mathbb{E}\left[\frac{1}{n} \sum_{i=1}^n||(\nabla f_{i}(u_{t-2})-\nabla f_{i} (x_*))||^2\right] + \left(1-\frac{1}{n}\right)\mathbb{E}\left[\frac{1}{n} \sum_{i=1}^n||(\nabla f_{i}(\phi_{i}^{t-2})-\nabla f_{i} (x_*))||^2\right] \notag \\
=& \cdots \notag \\
=& \sum_{j=1}^{t-1}\frac{1}{n}\left(1-\frac{1}{n}\right)^{t-1-j}\mathbb{E}\left[\frac{1}{n} \sum_{i=1}^n||(\nabla f_{i}(u_{j-1})-\nabla f_{i} (x_*))||^2\right] \notag \\
&+ \left(1-\frac{1}{n}\right)^{t-1}\left[\frac{1}{n} \sum_{i=1}^n||(\nabla f_{i}(\phi_{i}^0)-\nabla f_{i} (x_*))||^2\right]. \notag
\end{align}
Here we defined $\sum_{j=1}^0 = 0$ for $t=1$. 

Using these inequalities and the definition of $\eta$, we have 
\begin{align}
&\sum_{t=1}^{m_s} \left[ \frac{t}{\eta t-L_{\mathrm{max}}}\mathbb{E}||g_t - \nabla F(u_{t-1})||^2 -  \frac{1}{L_{\mathrm{max}}}\mathbb{E}\left[\frac{1}{n}\sum_{i=1}^{n}||\nabla f_i(u_{t-1})-\nabla f_i(x_*)||^2\right] \right] \notag \\
=&\sum_{t=1}^{m_s} \left[ \frac{1}{2L_{\mathrm{max}}}\mathbb{E}\left[\frac{1}{n}\sum_{i=1}^n ||(\nabla f_{i}(\phi_{i}^{t-1})-\nabla f_{i} (x_*))||^2\right]-  \frac{1}{2L_{\mathrm{max}}}\mathbb{E}\left[\frac{1}{n}\sum_{i=1}^n ||\nabla f_{i}(u_{t-1})-\nabla f_{i} (x_*)||^2\right]  \right] \notag \\
=&\frac{1}{2L_{\mathrm{max}}}\sum_{t=1}^{m_s} \left[ \sum_{j=1}^{t-1}\frac{1}{n}\left(1-\frac{1}{n}\right)^{t-1-j}\mathbb{E}\left[\frac{1}{n} \sum_{i=1}^n||(\nabla f_{i}(u_{j-1})-\nabla f_{i} (x_*))||^2\right] \right] \notag \\ 
&+ \frac{1}{2L_{\mathrm{max}}}\sum_{t=1}^{m_s}\left(1-\frac{1}{n}\right)^{t-1}\left[\frac{1}{n} \sum_{i=1}^n||(\nabla f_{i}(\phi_{i}^0)-\nabla f_{i} (x_*))||^2\right] \notag \\ &- \frac{1}{2L_{\mathrm{max}}}\sum_{t=1}^{m_s}\mathbb{E}\left[\frac{1}{n}\sum_{i=1}^{n}||\nabla f_i(u_{t-1})-\nabla f_i(x_*)||^2\right] .  \notag 
\end{align}
Observe that
\begin{align}
&\sum_{t=1}^{m_s} \left[ \sum_{j=1}^{t-1}\frac{1}{n}\left(1-\frac{1}{n}\right)^{t-1-j}\mathbb{E}\left[\frac{1}{n} \sum_{i=1}^n||(\nabla f_{i}(u_{j-1})-\nabla f_{i} (x_*))||^2\right] \right] \notag \\
=&\sum_{t=2}^{m_s}\frac{1}{n}\left(1-\frac{1}{n}\right)^{t-2}\mathbb{E}\left[\frac{1}{n} \sum_{i=1}^n||(\nabla f_{i}(u_{0})-\nabla f_{i} (x_*))||^2\right] \notag \\
&+ \sum_{t=3}^{m_s}\frac{1}{n}\left(1-\frac{1}{n}\right)^{t-3}\mathbb{E}\left[\frac{1}{n} \sum_{i=1}^n||(\nabla f_{i}(u_{1})-\nabla f_{i} (x_*))||^2\right] \notag \\
&+\cdots \notag \\
&+ \sum_{t=m_{s}}^{m_s}\frac{1}{n}\left(1-\frac{1}{n}\right)^{t-m_s}\mathbb{E}\left[\frac{1}{n} \sum_{i=1}^n||(\nabla f_{i}(u_{m_s-2})-\nabla f_{i} (x_*))||^2\right] \notag \\
\leq&\sum_{t=1}^{m_s}\mathbb{E}\left[\frac{1}{n}\sum_{i=1}^{n}||\nabla f_i(u_{t-1})-\nabla f_i(x_*)||^2\right]. \notag 
\end{align}
By Lemma \ref{object bound lemma} and the definition of $\phi_i^0$, we get
\begin{align}
&\sum_{t=1}^{m_s}\left(1-\frac{1}{n}\right)^{t-1}\left[\frac{1}{n} \sum_{i=1}^n||(\nabla f_{i}(\phi_{i}^0)-\nabla f_{i} (x_*))||^2\right] \notag \\
=&2m_s L_{\mathrm{max}}(P(x_0)-P(x_*)-\frac{\mu}{2}||x_0-x_*||^2).\notag  
\end{align}
Hence we get
\begin{align}
&\sum_{t=1}^{m_s} \left[\frac{t}{\eta t-L_{\mathrm{max}}}\mathbb{E}||g_t - \nabla F(u_{t-1})||^2 -  \frac{1}{L_{\mathrm{max}}}\mathbb{E}\left[\frac{1}{n}\sum_{i=1}^{n}||\nabla f_i(u_{t-1})-\nabla f_i(x_*)||^2\right] \right] \notag \\
\leq&m_s \left(P(x_0)-P(x_*)-\frac{\mu}{2}||x_0-x_*||^2\right).\notag
\end{align}
Combining Lemma \ref{main lemma2} with this result yields
\begin{align}
&\mathbb{E} [P({x}_{m_s})-P(x_{*})] + \frac{\eta+m_s\mu}{2m_s} \mathbb{E}||{v}_{m_s} - x_*||^2 \notag \\
\leq& \frac{1}{2}\left(P(x_0)-P(x_*)-\frac{\mu}{2}||x_0-x_*||^2\right)+ \frac{\eta}{2m_s}||{v}_{0} - x_*||^2. \notag
\end{align}
Since $x_{m_s}=\widetilde{x}_s$, $v_{m_s}=\widetilde{v}_s$, $x_0=\widetilde{x}_{s-1}$, and  $v_0=(1-\alpha)\widetilde{v}_{s-1}+\alpha\widetilde{x}_{s-1}$, we obtain
\begin{align}
&\mathbb{E} [P(\widetilde{x}_{s})-P(x_{*})] + \frac{\eta+m_s\mu}{2m_s} \mathbb{E}||\widetilde{v}_{s} - x_*||^2 \notag \\
\leq& \frac{1}{2}(P(\widetilde{x}_{s-1})-P(x_*)) + \left(\frac{\alpha \eta}{2m_s}-\frac{\mu}{4}\right)||\widetilde{x}_{s-1} - x_*||^2 \notag \\ &+\frac{(1-\alpha)\eta}{2m_s}||\widetilde{v}_{s-1} - x_*||^2. \notag
\end{align}
Finally, taking expectations with respect to all previous stages gives the desired result. 
\qed
\end{pr main thm2}
\end{document}